\documentclass[11pt]{article} 

\usepackage[colorlinks=true,linkcolor=blue,urlcolor=blue]{hyperref}
\usepackage{dsfont}

\usepackage{amsfonts,latexsym,amstext} 
\usepackage{amsmath} 
\usepackage{amssymb} 
\usepackage{comment} 
\usepackage{amsthm} 
\usepackage[latin1]{inputenc} 
\usepackage{color} 
 
\theoremstyle{plain} 
\newtheorem{theorem}{Theorem}[section] 
\newtheorem{lemma}[theorem]{Lemma} 
\newtheorem{proposition}[theorem]{Proposition} 
\newtheorem{corollary}[theorem]{Corollary} 

\newtheorem{definition}[theorem]{Definition} 
\newtheorem{hypothesis}[theorem]{Hypothesis} 
\theoremstyle{remark} 
\newtheorem{remark}[theorem]{Remark} 

\allowdisplaybreaks 
 
\makeatletter 
\@addtoreset{equation}{section} 
 
\makeatother 
 
\def\sqr#1#2{{\vcenter{\vbox{\hrule height .#2pt \hbox{\vrule 
 width .#2pt height#1pt \kern#1pt \vrule 
width .#2pt} \hrule height .#2pt}}}}

\def\ds{\begin{displaystyle}} 
\def\eds{\end{displaystyle}} 
 
\def\<{\langle } 
\def\>{\rangle }

\def\R{\mathbb R} 
\def\N{\mathbb N}

\def\E{\mathbb E} 
\def\P{\mathbb P} 
\def\Q{\mathbb Q}

\newcommand{\sper}[1]{\mathbb{E} \left[ #1 \right]}                               
 
\DeclareMathAlphabet{\mathonebb}{U}{bbold}{m}{n}                           %
\newcommand{\one}{\ensuremath{\mathonebb{1}}}                               
 
\title{The identification problem for BSDEs driven by possibly non quasi-left-continuous random measures} 

\author{ 
Elena Bandini 
\thanks{Universit\`a degli Studi di Milano Bicocca, Dipartimento di Matematica, via  Roberto Cozzi 55, 20125 Milano, Italy; e-mail: elena.bandini@unimib.it.} 
\and Francesco Russo 
\thanks{Unit\'e de Math\'ematiques appliqu\'ees, ENSTA ParisTech, Universit\'e Paris-Saclay, 
828, boulevard des Mar\'echaux, F-91120 Palaiseau, France; e-mail: francesco.russo@ensta-paristech.fr.} 
} 
\date{} 
 
\begin{document} 
 
\allowdisplaybreaks 
\maketitle

 \begin{abstract} 
In this paper we focus on the so called \emph{identification problem} for a backward SDE driven by a continuous  local martingale
and a possibly non quasi-left-continuous random measure.
Supposing that a solution $(Y,Z,  U)$ of a backward SDE
is such that $Y_t = v(t,X_t)$ where $X$ is  an underlying process and 
$v$ is a deterministic function, solving the identification problem consists
in determining $Z$ and $U$ in term of $v$. 
We study the over-mentioned identification problem
under various sets of assumptions and we provide a family of examples
including the case when $X$ is a non-semimartingale jump process
solution of an SDE with singular coefficients.
\end{abstract}

 
 \medskip 
 

\section{Introduction}

 This paper considers a  BSDE  driven by a compensated random 
measure $\mu-\nu$, of the form 
\begin{align}\label{GeneralBSDE} 
Y_t &= \xi 
+\int_{]t,\,T]\times \R} \tilde f(s,\,e,\,Y_{s-},\,Z_s,\,U_{s}(e))\,  d \zeta_s\nonumber\\ 
&\quad- \int_{]t,\,T]} Z_s \, dM_s - \int_{]t,\,T]\times \R} U_s(e)\,(\mu-\nu)(ds\,de), 
\end{align} 
whose solution is a triplet of processes   $(Y,Z,U)$, with $Y$ a c\`adl\`ag adapted process, $Z$ a predictable process and $U(\cdot)$ a predictable random field. 
Besides $\mu$ and $\nu$ appear two driving random elements, namely a continuous martingale $M$ and 
a non-decreasing adapted c\`adl\`ag process $\zeta$,
while $\xi$ is a  square integrable random variable, and 
$ \tilde f$ is a  random function. 
 Often $Y$ turns out to be 
of the type $v(t,X_t)$ where $v$ is a deterministic function, and $X$ is  a
 c\`adl\`ag adapted process.  
The \emph{identification problem} consists in determining $Z$ and $U$ in terms of $v$.

BSDEs  have been deeply studied   since the seminal paper 
\cite{pardouxpeng}, where the  Brownian context  appears as a particular case of  \eqref{GeneralBSDE},  setting 
$\mu
=0$, $\zeta_s \equiv s$. 
There, $M$ is a standard Brownian 
motion and $\xi$ is measurable with respect to the Brownian $\sigma$-field 
at terminal time. 
 In that case  the unknown can be reduced to $ (Y,Z)$, since   $U$ can be arbitrarily chosen. 
BSDEs with a discontinuous driving term of the form \eqref{GeneralBSDE}  have been studied as well; in almost all cases, the random measure $\mu$  is quasi-left-continuous, i.e.
$\mu(\{S \} \times \R)=0$ on $\{S < \infty \}$ for every predictable time $S$, 
 see, e.g.,  \cite{xia}, \cite{BP}, \cite{TaLi}, \cite{BandiniFuhrman}, 
  \cite{BandiniConfortola}. 
Existence and uniqueness for 
BSDEs driven  by a 
random measure 
 which is not necessarily quasi-left-continuous 
are very recent, 
and were 
 discussed   in  \cite{BandiniBSDE} in the purely discontinuous case, and in \cite{papapantoleon_possamai_saplaouras} in the jump-diffusion case. 
 
When the random dependence 
of   $\tilde f$ 
is provided by a Markov solution $X$ of a 
 forward  SDE, and $\xi$ is a real function of $X$ at the terminal time $T$, then the   BSDE \eqref{GeneralBSDE} 
 is called forward BSDE. 
 In the Brownian context, 
when $X$ is the solution of a classical SDE with 
diffusion coefficient $\sigma$, forward BSDEs  generally constitute  stochastic representations of a 
partial differential equation. If  $v:[0,T] \times \R \rightarrow \R$ is a 
classical (smooth) solution of the mentioned  
PDE, 
then $Y_s = v(s,X_s)$, $Z_s = \sigma(s,X_s)\, \partial_x v(s,X_s)$, 
 generate  a solution 
to the forward BSDE, 
see e.g. \cite{Peng91}, \cite{PardouxPeng92}, \cite{Peng92b}, which provide the solution to  the identification problem in that particular case. 
Conversely, 
solutions of forward BSDEs generate solutions  of PDEs in the viscosity sense,
or in other generalized sense, see e.g.  \cite{barrasso3}. 
More precisely, for each given couple $(t,x) \in [0,\,T]\times \R$, consider an  underlying process $X$ given by 
the solution $X^{t,x}$ of an SDE starting at $x$ at time $t$;
if $(Y^{t,x},Z^{t,x})$ is a family of solutions of the forward BSDE,  
 under reasonable general assumptions, the function $v(t,x):= Y_t^{t,x}$ is 
a viscosity solution of the related PDE. 
In the Brownian context,  the identification problem of $Z$ has been faced even if $v \in C^{0,1}$,   including  the infinite dimensional case, see 
for instance  \cite{FuhrmanTessitore} under slightly more general conditions.

In the general case when 
the forward  BSDEs are 
also  driven by random measures,  these equations  generally constitute  stochastic representations of a 
partial integro-differential equation (PIDE). 
When $v$ is a classical solution of the PIDE and $X$ is a solution to a Markov  jump-diffusion equation,    the identification problem was solved  in  \cite{BaBuPa}. 
Analogous results can be obtained when  $X$ is a  purely discontinuous Markov process,  see \cite{CoFu-m}.
In both cases,  the BSDE is  driven   by a compensated random measure $\mu-\nu$ with  $\mu$  quasi-left-continuous.
In the context of alternative BSDEs with jumps, i.e. the one	
of martingale driven forward BSDEs of \cite{CarboneFerrSantacroce},
the identification problem was discussed for instance in
 \cite{LaachirRusso} and \cite{barrassopreprint1}.

In \cite{BandiniRusso2} we  extended the above-mentioned identification results in two directions. 
Firstly,  we generalized
 \cite{CoFu-m} to the case of  BSDEs driven by non-quasi-left-continuous random measures, related to a special class of piecewise deterministic Markov process (PDMPs).
Secondly, in the  non purely discontinuous  case, we  extended the study  to the case when $Y_t = v(t, X_t)$, with $X$  a special weak Dirichlet process of finite quadratic variation and $v$  of class $C^{0,1}$.
A similar technique was used earlier in the different context of 
verification theorems for control problems, see \cite{gr1}.

Besides the survey aspects, the present paper extends the results of \cite{BandiniRusso2} along three lines.
\begin{enumerate}
\item One investigates the identification problems, going beyond the forward BSDEs formalism, even though following the same lines of \cite{BandiniRusso2}. 
\item We generalize the results in \cite{BandiniRusso2} by alleviating some important hypotheses. 
As a matter of fact, Proposition  \ref{L_ident_mu_muX}   improves the achievements of Proposition 2.17  in  \cite{BandiniRusso2}, since the condition 
\begin{equation}\label{CondInMeno}
\nu(\{S\},de)=\mu(\{S\},de)\,\textup{a.s.} \,\,\textup{for every predictable time}\,\, S \,\,\textup{such that}\,\, [[S]] \subset K
\end{equation}
	is no longer needed here. 
	 This allows to formulate Theorems \ref{P_ident} and \ref{P_ident_C0} under more general assumptions, and 
	 extends the applicability of our results.  Among others, 
	 we are  able to solve the identification result for more general jump-diffusion processes and   piecewise deterministic Markov processes, see respectively
	 Corollaries  \ref{C_id_BBP}  and \ref{C_id_PDPs}. 
\item 
We apply our results  to the completely new case when $X$ is the solution to a martingale problem with general jumps and distributional drift,  related thus
to an operator of the form $\beta'(x)\frac{\partial}{\partial x} + \frac{1}{2}\sigma^2(x) \frac{\partial^2}{\partial x^2}$, with $\beta$ only continuous, and to some predictable random measure $\nu$  possibly discontinuous.  Martingale problems of this type have been studied in the companion paper \cite{BandiniRusso_DistributionalDrift}.  
We  solve the identification problem in this context, see  Corollary \ref{C_id_distr}.


\end{enumerate}

\section{Notation and preliminaries}\label{S:2}
We fix a   positive horizon $T$. Given a topological space $E$, in the sequel $\mathcal{B}(E)$ will denote 
the Borel $\sigma$-field associated with $E$. We will indicate by 
 $C^{0,1}$ 
the space of all functions 
$u: [0,\,T]\times \R \rightarrow \R$,  $(t,x)\mapsto u(t,x)$, 
that are continuous  together their derivative 
$\partial_x u$.

For a filtered probability space $(\Omega,\mathcal{F}, (\mathcal{F}_t)_{t \geq 0},\P)$, we will always suppose that
$(\mathcal{F}_t)_{t \geq 0}$ satisfies the usual conditions, with $\mathcal F = \mathcal F_T$. 
Related to it, 
$\mathcal{P}$ (resp. $\mathcal{\tilde{P}}=\mathcal{P}\otimes \mathcal{B}(\R)$) will denote the predictable $\sigma$-field on $\Omega \times [0,T]$ (resp. on $\tilde{\Omega}= \Omega \times [0,T]\times \R$). 
Analogously, we set $\mathcal{O}$ (resp. $\tilde{\mathcal{O}} =\mathcal{O} \otimes \mathcal{B}(\R)$) as the optional   $\sigma$-field on $\Omega \times [0,T]$ (resp. on $\tilde{\Omega}$). 
Moreover, 
$\tilde {\mathcal{F}}$ will be $\sigma$-field $\mathcal{F} \otimes {\mathcal B}([0,T] \times \R)$, and 
we will  indicate by $\mathcal{F^{\P}}$ 
the completion of $\mathcal F$ with the $\P$-null sets. 
We set 
$\mathcal{\tilde F}^{\P} =  \mathcal{F}^{\P} \otimes {\mathcal B}([0,T] \times \R)$. 
By default, all the stochastic processes will be considered with parameter $t \in [0,\,T]$. 
By convention, any c\`adl\`ag process defined on $[0,\,T]$ is  extended to $\R_+$ by continuity.
A random set $A \subset \tilde \Omega$ is called evanescent if the set $\{\omega: \exists \,t \in \R_+ \textup{ with } (\omega, t) \in A\}$ is $\P$-null. Generically, all the equalities of random sets will be intended up to an evanescent set. 

For a measurable process $X$ we denote by ${}^p{} 
(X)$ its predictable projection, see e.g. Theorem 5.2 in \cite{chineseBook}.
 A bounded variation process $X$ on $[0,\,T]$  will be said to be with integrable variation if the expectation of its total variation is finite. 
$\mathcal{A}$ (resp. $\mathcal{A}_{\textup{loc}}$) will denote  the collection of all adapted processes with   integrable variation (resp.  with locally integrable variation), and    $\mathcal{A}^+$ (resp $\mathcal{A}_{\textup{loc}}^+$)  the collection of all adapted integrable increasing (resp. adapted locally integrable)  processes. 
In general, these notions refer to the underlying probability $\P$; when this is not the case, we will mention the specific probability. 
The significance of locally is the usual one which refers 
to  localization by stopping times, see e.g. (0.39) of 
\cite{jacod_book}. 

 We also recall that a random kernel $\phi(a, db)$ of a measurable space $(A, \mathcal A)$ into another measurable space $(B, \mathcal B)$ is a family $\{\phi(a, \cdot), \,\,a \in A\}$ of positive measures on $(B, \mathcal B)$, such that  $\phi(\cdot, C)$ is $\mathcal A$-measurable for any $C \in \mathcal B$. 
Finally, the concept of random measure  
is extensively used 	throughout   the paper:  
for a  detailed  discussion on this topic  and the unexplained  notations see 
Chapter I and Chapter II, Section 1, in \cite{JacodBook}, Chapter III in \cite{jacod_book},  and  Chapter XI, Section 1, in \cite{chineseBook}.
In particular, if $\mu$ is a random measure on $[0,\,T]\times \R$, for any measurable real function $H$ defined on $\tilde \Omega$, one denotes $H \star \mu_t:= \int_{]0,\,t] \times \R} H(\cdot, s,x) \,\mu(\cdot, ds \,de)$.

\subsection{Stochastic integration with respect to integer-valued random measures} \label{Sec_2.1}
Let  $(\Omega,\mathcal{F}, (\mathcal{F}_t)_{t \geq 0},\P)$ be a filtered probability space.
In the  sequel of  Section \ref{S:2}, $\mu$ will be  an integer-valued random measure on $[0,\,T] \times \R$, and $\nu$ will be  its compensator, for which we will  choose the "good" version 
as constructed in point (c) of 
Proposition 1.17, Chapter II, in \cite{JacodBook}.
Set 
\begin{align}
D&=\{(\omega, t): \mu(\omega, \{t\}\times \R)>0\}, \label{D}\\ 
J&=\{(\omega, t): \nu(\omega, \{t\}\times \R)>0\},\label{J}\\
K&=\{(\omega, t): \nu(\omega, \{t\}\times \R)=1\}.\label{K}
\end{align}
We define $\nu^d:= \nu\,\one_{J}$ and $\nu^c:= \nu\,\one_{J^c}$. 
Similar conventions will be used for the other integer-valued random measures appearing in the paper.
We will sometimes  use the form of $\mu$ given  in Proposition 1.14, Chapter II,  in \cite{JacodBook}, i.e.
\begin{equation}\label{mubeta}
\mu(dt\,de) = \sum_{s \geq 0}\one_{D}(s,\omega)\,\delta_{(s,\beta_s(\omega))}(dt\,de), 
\end{equation}
where $\beta$ is a real-valued optional process.
\begin{remark}\label{R_pred_supp}
\begin{itemize}
\item[(i)]	$D$ is a thin set, namely $D=\cup_n[[T_n]]$ with $(T_n)_n$ random times, see Theorem 11.13 in \cite{chineseBook}.
\item [(ii)]	
	$J$ is the predictable support of $D$, namely $J = \{{}^p{}(\one_D)>0\}$, see  Theorem 5.39 in \cite{chineseBook}. 
This is equivalent to  
$\one_J= 
{}^p{}(\one_D)$.
\item[(iii)]
 There exists a sequence of predictable
times $(R_n)_n$ with disjoint graphs, such that $J =
\cup_n[[R_n]]$, see   Proposition
2.23, Chapter I, in \cite{jacod_book}.
 \item[(iv)]
	$K$ is the largest predictable subset of $D$,
	 see  
	 Theorems 11.14 in \cite{chineseBook}. 
Since $K$ is predictable, 
we have 
${}^p{} 
(\one_K ) 
= 
\one_K$. 
\item[(v)]	A progressive set $B$ contained in a thin set is also  thin set,   	
	 see Theorem 3.19 in \cite{chineseBook}. In particular, $K$ is a thin set.

\end{itemize}

\end{remark} 

\begin{remark}\label{R_H_mu} 
$\nu$ admits a disintegration of the type 
		\begin{equation}\label{nu_dis} 
		\nu(\omega, ds\,de)=d A_s(\omega) \,\phi(\omega, s,de), 
		\end{equation} 
		where  $\phi$ is a random kernel from $(\Omega \times [0,\,T], \mathcal P)$ into $(\R, \mathcal B(\R))$ and  $A$ is a right-continuous nondecreasing predictable process, such that $A_0=0$, see 
		for instance Remark 4.4 in \cite{BandiniBSDE}. 
\end{remark} 

We recall an important notion of measure associated with  $\mu$, given in formula (3.10) in \cite{jacod_book}. 
\begin{definition} 
	\label{D_DoleansMeas} 
	Let $(\tilde \Omega_n)$ be a  partition  of $\tilde \Omega$ constituted by elements of $\tilde{\mathcal O}$, such that $\one_{\tilde \Omega_n} \star \mu \in \mathcal A$. 
	$M^{\P}_{\mu}$ denotes the $\sigma$-finite measure on  $(\tilde{\Omega}, 
	\mathcal{{\tilde F}^{\P}})$, such that 
	for every $W: \tilde \Omega \rightarrow \R $ positive, bounded,   $\mathcal{\tilde{F}}^{\P}$-measurable function, 
	\begin{equation}\label{DoleansMeas} 
	M^{\P}_{\mu}(W\,\one_{\tilde{\Omega}_n})= \E\big[W\,\one_{\tilde{\Omega}_n} \star \mu_T\big]. 
	\end{equation} 	 
\end{definition} 
Let us  set $\hat{\nu}_t(de):= \nu(\{t\}, de)$ for all $t \in [0,\,T]$. 
For 
any 
$W \in \mathcal{\tilde{O}}$, 
we define 
\begin{align*} 
	\hat{W}_t = \int_{\R} W_t(e)\,\hat{\nu}_t(de), 
	\quad \tilde{W}_t = \int_{\R} W_t(e)\,\mu(\{t\}, de)- \hat{W}_t,\quad t \geq 0,\nonumber
\end{align*}
with the convention that
	$\tilde{W}_t= + \infty$ if $\hat{W}_t$ is not defined.  
For every $q \in [1,\,\infty[$, we  introduce the linear spaces 
\begin{align*}
	\mathcal{G}^q(\mu)=\Big\{W \in \mathcal{\tilde{P}} 
	:\,\, 
	\Big[\sum_{s\leq \cdot}|\tilde{W}_s|^2\Big]^{q/2}\in \mathcal{A}^+\Big\},\quad 
	\mathcal{G}^q_{\textup{loc}}(\mu)=\Big\{W \in  \mathcal{\tilde{P}} 
	: \,\, 
	\Big[\sum_{s\leq \cdot}|\tilde{W}_s|^2\Big]^{q/2}\in \mathcal{A}_{\textup{loc}}^+\Big\}.\nonumber
\end{align*} 
Given $W \in \mathcal{\tilde{P}}$, we define 
the  increasing (possibly infinite) predictable process
\begin{align}
C(W) := |W - \hat{W}\,\one_{J}|^2\star \nu + \sum_{s \leq \cdot}(1-\hat{\nu}_s(\R))\,|\hat{W}_s|^2\,\one_{J \setminus K}(s), \label{C(W)BIS}
\end{align}
provided the right-hand side is well-defined.
By  
Theorem 11.21, point 3)   in \cite{chineseBook}, if $W \in \mathcal{G}^2(\mu)$, then $\langle W \star (\mu-\nu)\rangle$ is well defined and  
\begin{align}
	C(W)= \langle W \star (\mu-\nu)\rangle.
	\label{bracket_mart} 
\end{align}
We set
\begin{align}\label{G2norm}
||W||^2_{\mathcal{G}^2(\mu)}:=\sper{C(W)_T}.	
\end{align}
We also introduce the  space 
\begin{equation*}\label{L2mu} 
	\mathcal{L}^2(\mu):=\left\{W \in \tilde{\mathcal{P}} 
	:\,\,  ||W||_{\mathcal{L}^2(\mu)}:=\E\Big[\int_{]0,T]\times \R}  |W_s(e)|^2 \,\nu(ds\,de)\Big] < \infty\right\}.
\end{equation*}

\begin{remark}[Lemma 2.4 in \cite{BandiniRusso2}]\label{L_G2mu_L2mu_inclusion} 
	If $W \in \mathcal{L}^2(\mu)$ 
then $W \in \mathcal{G}^2(\mu)$ 
	 and 
	$||W||^2_{\mathcal{G}^2(\mu)} \leq ||W||^2_{\mathcal{L}^2(\mu)}$.
Moreover   $\mathcal{L}^2_{\rm loc}(\mu) \subset \mathcal{G}_{\rm loc}^2(\mu)$. 	 
\end{remark}

The following result is the object of Proposition 2.8 in \cite{BandiniRusso2}.
 \begin{proposition}\label{P_forBSDEs2} 
If 
$C(W)_T=0$ a.s.,
then 
$||W-\hat{W}\,\one_{K}||_{\mathcal{L}^2(\mu)}=0$,
or, equivalently,  there exists a predictable process $(l_s)$ such that
\begin{equation}\label{2.12}
W_s(e)= l_s\,\one_{K}(s),\quad d\P\,\nu(ds\,de)\textup{-a.e.}
\end{equation}
In particular, 
$W_s(e)=0$, $d\P\,\nu^c(ds\,de)$\textup{-a.e.}, 
and  there is a predictable process  $(l_s)$  such that \\
$W_s(e)= l_s\,\one_{K}(s)
$, $d\P\,\nu^d(ds\,de)$\textup{-a.e.}	
\end{proposition} 
We recall the following definition, fundamental for the sequel of the paper. 
\begin{definition} Given a c\`adl\`ag  process $X$, we introduce the associated  jump measure,  namely 
the integer-valued random measure on $\R_+ \times \R$ defined as 
\begin{equation}\label{X_jump_measure} 
	\mu^X(\omega; dt\,dx):= \sum_{s\in ]0,\,T]} \,\one_{\{\Delta X_s(\omega)\neq 0\}}\,\delta_{(s,\,\Delta X_s(\omega))}(dt\,dx). 
\end{equation} 
The compensator of 	$\mu^X(ds\,dx)$ 
will be denoted by $\nu^X(ds\,dx)$.
\end{definition}
We will consider the following condition for a couple $(\chi, \Q)$ where $\chi$ is a random measure and $\Q$ is a given probability (of simply for $\chi$ when $\Q$ is the self-explanatory probability $\P$):
\begin{equation}\label{Sec:WD_CNS} 
\int_{]0,\cdot]\times \R} (|x|\wedge |x|^2)\,\chi(ds\,dx) \in \mathcal{A}^{+}_{\rm loc}. 
\end{equation} 
A more restrictive condition will be also considered, namely 
\begin{equation}\label{alpha+1integ} 
\int_{]0,\cdot]\times \R} (|x|\wedge |x|^{1+\alpha})\,\chi(ds\,dx) \in \mathcal{A}^{+}_{\rm loc}\quad \textup{for some}\,\, \alpha \in [0,\,1].
\end{equation} 
By abuse of notations, when $\chi = \mu^X$ for a given  c\`adl\`ag   process $X$, we will say that $X$ verifies condition \eqref{Sec:WD_CNS}  or condition \eqref{alpha+1integ} (under $\Q$).
\begin{remark}\label{R:WD_CNS}
Condition \eqref{alpha+1integ} 
implies  in particular 
condition \eqref{Sec:WD_CNS}. 
\end{remark}

We will be also interested in functions  $v:[0,T] \times \R \rightarrow \R$  of class $C^{0,1}$ fulfilling the  integrability property (with respect to $(\chi,\Q)$ or simply for $\chi$ when $\P$ is self-explanatory)
\begin{equation}\label{Sec:WD_E_C} 
\int_{]0,\cdot]\times \R} |v(s,X_{s-} +x )-v(s,X_{s-})-x\,\partial_x v(s,X_{s-})|\,\one_{\{|x| >1\}}\,\chi(ds\,dx) \in \mathcal{A}^{+}_{\rm loc}. 
\end{equation}
Also in this case, by abuse of notations, when $\chi = \mu^X$ for a given  c\`adl\`ag   process $X$, we will say that $v$ and $X$ verify condition \eqref{Sec:WD_E_C}   (under $\Q$).

\subsection{Chain rules for special weak Dirichlet processes }

Let  $(\Omega,\mathcal{F}, (\mathcal{F}_t)_{t \geq 0},\P)$ be a filtered probability space. Special weak Dirichlet processes constitute a further development of weak Dirichlet processes, which were introduced by \cite{er2}, \cite{gr}
in the continuous case and by \cite{cjms}
in the jump case.

 $O$  is  an $({\mathcal F}_t)$-orthogonal process if 
	$[O, N] = 0$  for every $N$ continuous local $({\mathcal F}_t)$-martingale. We recall that $[\cdot, \cdot]$ is the covariation extending the classical covariation of semimartingales, see \cite{rv95} and Definition 2.4 in \cite{BandiniRusso1}. An $({\mathcal F}_t)$-local  martingale $M$ is said to be purely discontinuous if  it is $({\mathcal F}_t)$-orthogonal. 	
A special weak Dirichlet process is a process of the type $X= M+ A$, where $M$ is an $(\mathcal F_t)$-local martingale and $A$ is an $(\mathcal F_t)$-predictable orthogonal process, see Definition 
5.6 in \cite{BandiniRusso1}. 
When $A$ has bounded variation, then $X$ is a special $(\mathcal F_t)$-semimartingale. 
	\begin{remark}[Proposition 5.9 in \cite{BandiniRusso1}]\label{Sec:WD_P_unique_decomp}
	Any $(\mathcal F_t)$-special weak Dirichlet process $X$ admits a unique 
 decomposition of the type
 \begin{equation}\label{Mc+Md+A} 
	X=X^c+M^d+A, 
	\end{equation} 
	where $X^c$ is a continuous local martingale, $M^d$ is a purely discontinuous local martingale, and $A$ is an $({\mathcal F}_t)$-predictable and orthogonal process, with $A_0=0$.
	\eqref{Mc+Md+A} is called the canonical decomposition of $X$.
\end{remark}

In the sequel we will   consider the following  assumptions on a couple $(X,Y)$ of adapted processes.
\begin{hypothesis}\label{H_chain_rule_C01} 
	$X$ is an $(\mathcal F_t)$-special weak Dirichlet process      of finite quadratic variation with its canonical decomposition  $X=X^c+M^d+A$, satisfying condition \eqref{Sec:WD_CNS}. 
	$Y_t = v(t,\,X_t)$ for some (deterministic) function $v: [0,T] \times \R \rightarrow \R$ of class $C^{0,1}$ such that $v$ and $X$ verify   condition \eqref{Sec:WD_E_C}. 
\end{hypothesis} 
\begin{remark}
\begin{itemize}
\item[(i)] [Proposition 4.5 in \cite{BandiniRusso1}]. If $X$ is a c\`adl\`ag  process such that $\sum_{s \leq T}|\Delta X_{s}|^2 < \infty $ a.s.,
	 and     $v:[0,T] \times \R \rightarrow \R$  is  a function of class $C^{0,1}$, then	 
	\begin{align}
		|v(s,X_{s-} + x)-v(s,X_{s-})|^2\,\one_{\{\vert x \vert \leq 1\}}\star \mu^X \in \mathcal{A}^+_{\rm loc}.\label{A2_Aloc_NEW} 
		\end{align}
\item[(ii)] [Lemma 5.29 in \cite{BandiniRusso1}]. 
If $X$ is a c\`adl\`ag  process  satisfying 
	condition \eqref{Sec:WD_CNS}, and     $v:[0,T] \times \R \rightarrow \R$  is  a function of class $C^{0,1}$ fulfilling \eqref{Sec:WD_E_C}, then	 
\begin{align}
		|v(s,X_{s-} + x)-v(s,X_{s-})|\,\one_{\{\vert x \vert > 1\}}\star \mu^X \in \mathcal{A}^+_{\rm loc}.\label{F2_Aloc_NEW} 
		\end{align}	
\item[(iii)][Remark 5.30 in \cite{BandiniRusso1}]. Condition \eqref{Sec:WD_E_C} is automatically verified if $X$ is a c\`adl\`ag process satisfying
\eqref{Sec:WD_CNS}  and $v:[0,T] \times \R \rightarrow \R$  is  a function of class $C^{0,1}$ with $\partial_x v$ bounded.
\end{itemize} 
\end{remark}

\begin{theorem}[Theorem 5.31  in \cite{BandiniRusso1}]\label{T_C1_dec_specialweak_Dir} 
Let $(X,Y)$ be a couple of $(\mathcal F_t)$-adapted processes satisfying Hypothesis  \ref{H_chain_rule_C01}  with corresponding function $v$.
		Then  
		we have 
	\begin{align}\label{C01_special_WD_formula} 
	v(t,X_t)&=v(0,X_0)+\int_0^t\partial_x v(s,X_s)\,d X^c_s\nonumber\\ 
	&+ \int_{]0,\,t]\times \R} (v(s,X_{s-}+x)-v(s,X_{s-}))\,(\mu^X-\nu^X)(ds\,dx)+A^v(t), 
	\end{align} 
	where 
	$A^v$ is a predictable $(\mathcal F_t)$-orthogonal  process. 
\end{theorem}

We now need to formulate a technical assumption under which Proposition \ref{P_C00_chain_rule} below  holds, see item (iii) of Hypothesis \ref{H_chain_rule_C0}. Let $E$ be a closed subset of $\R$ on which $X$ takes values. Given a  c\`adl\`ag function $\varphi: [0,\,T] \rightarrow \R$, we denote by $\mathcal C_{\varphi}$ the set of times $t \in [0,\,T]$ for which there is a left (resp. right) neighborhood $I_{t-} = ]t-\varepsilon,\,t[$ (resp. $I_{t+} = [t,\,t+\varepsilon[$) such that $\varphi$ is  constant on $I_{t-}$ and $I_{t+}$.
Let us then consider the following    assumptions on a couple $(X,Y)$ of adapted processes.
\begin{hypothesis}\label{H_chain_rule_C0} 
	\hspace{2em} 
	\begin{itemize} 
		\item[(i)] 
		$Y$ is an $(\mathcal F_t)$-orthogonal process such that $\sum_{s \leq T}|\Delta Y_s| < \infty$, a.s.
		\item[(ii)]$X$ is a c\`adl\`ag process and  $Y_t = v(t,\,X_t)$ for 
		some 
		deterministic function $v: [0,T] \times \R \rightarrow \R$, satisfying  the  integrability condition 
		\begin{equation}\label{EC_C0} 
		\int_{]0,\,\cdot]\times \R}|v(t,X_{t-}+x)-v(t,X_{t-})|\,\mu^X(dt\,dx)\,\,\in \mathcal{A}^+_{\rm loc}. 
		\end{equation} 
		\item[(iii)] There exists $\mathcal C \in [0,\,T]$ such that for $\omega$ a.s. $\mathcal C \supset \mathcal C_X(\omega)$, and 
		\begin{itemize}
		\item[(i)] $\forall t \in \mathcal C$, $t \mapsto v(t,x)$ is continuous for all $x \in E$;
		\item[(ii)] $\forall t \in \mathcal C^c$, $x \in E$,  $(t, x)$ is a continuity point of $v$.
		\end{itemize}
	\end{itemize} 
\end{hypothesis}
\begin{remark}\label{R_1BIS}
	Item (iii) of Hypothesis \ref{H_chain_rule_C0}  is fulfilled in  two typical situations.
	\begin{enumerate}
	\item$\mathcal C=[0,\,T]$. Almost surely $X$ admits a finite number of jumps and $t \mapsto v(t,x)$ is continuous for all $x \in E$.
	\item$\mathcal C=\emptyset$ and  $v|_{[0,\,T] \times E}$ is continuous.
	\end{enumerate}
\end{remark}

\begin{proposition}[Proposition 5.37 
in \cite{BandiniRusso1}]\label{P_C00_chain_rule} 
Let $(X,Y)$ be a couple of $(\mathcal F_t)$-adapted processes satisfying Hypothesis  \ref{H_chain_rule_C0} with corresponding function $v$.
Then $v(t,X_t)$ is an $(\mathcal F_t)$-special weak Dirichlet c\`adl\`ag process with decomposition 
\begin{equation}\label{C00_decomp} 
v(t,X_t)=v(0,X_0) 
+ \int_{]0,\,t]\times \R} (v(s,X_{s-}+x)-v(s,X_{s-}))\,(\mu^X-\nu^X)(ds\,dx)+A^v(t), 
\end{equation} 
where $A^v$ is a predictable $(\mathcal F_t)$-orthogonal  process. 
\end{proposition}

 \section{A class of stochastic processes $X$ related in a specific way to an integer-valued random measure $\mu$}\label{Sec_2.2}

Let  $(\Omega,\mathcal{F}, (\mathcal{F}_t)_{t \geq 0},\P)$ be a filtered probability space.
We will make use of the following assumption relating a  c\`adl\`ag process $X$ and an integer-valued random measure $\mu(ds\,de)$ on $[0,T]\times \R$, with compensator $\nu(ds\,de)$
 (the random sets $D$, $J$ and $K$ 
  are defined in \eqref{D}-\eqref{J}-\eqref{K}). 
\begin{hypothesis}\label{H_X_mu} 
	We suppose that $X$  is an adapted c\`adl\`ag process  with decomposition $X=X^i + X^p$, where  the conditions below hold. 
\begin{enumerate}
	\item[1.a)] $Y:=X^i$ is a c\`adl\`ag quasi-left-continuous adapted process satisfying	$\{\Delta Y \neq 0\} \subset D$. 
	\item[1.b)] There exists a $\tilde{\mathcal{P}}$-measurable map  $\tilde \gamma: \Omega \times ]0,\,T] \times \R\rightarrow \R$ such that 
	\begin{equation}\label{B35ii} 
	\Delta Y_t(\omega)\,\one_{]0,\,T]}(t)=   \tilde \gamma(\omega,t,\cdot)\quad dM^\P_{\mu}\textup{-a.e.} 
	\end{equation} 
  \item[2.] $X^p$ is a c\`adl\`ag predictable process satisfying 	$\{\Delta X^p \neq 0\} \subset J$. 
\end{enumerate}
\end{hypothesis}
\begin{remark}\label{R:totinac_leftcont}
	We recall that  a random time $T$ is totally inaccessible if 
	$
	\one_{[[T]]}(\omega, S(\omega))\,\one_{\{S< \infty\}}= 0
	$ 
	for every predictable random time $S$
	see  Definition 2.20, Chapter I,  in \cite{JacodBook}, while a process $Y$ is quasi-left-continuous if $\Delta Y_S \one_{S < \infty}=0$ for all predictable random
time $S$, see   Definition 2.25, Chapter I, in \cite{jacod_book}.
	$Y$ is quasi-left-continuous	 if and only if there is a sequence of totally inaccessible times  $(T_n)$, with $[[T_n]] \cap [[T_m]] = \emptyset$, $n \neq m$, such that $\{\Delta Y \neq 0 \} = \cup_n [[T_n]]$, 
see Proposition 2.26, Chapter I,  
in \cite{JacodBook}. 
\end{remark}

In the sequel we will also need  the following assumption on  $\mu$.

\begin{hypothesis}\label{H_nu} 
		$J=K$ (up to an evanescent set).
\end{hypothesis} 
\begin{remark}\label{R_Hp_3.1} 
	Hypothesis \ref{H_nu},   is equivalent to ask that 
	\begin{align}\label{old_H}
D\,\,\textup{ is the disjoint union of }\,\, K \,\, \textup{and}\,\,  \cup_n [[T^i_n]] \,\, \textup{(up to an evanescent set)}\notag\\
\textup{with}\,\, (T^i_n)_n\,\, \textup{disjoint totally inaccessible times.}
\end{align}		
		 Indeed, if  \eqref{old_H} holds, then Hypothesis \ref{H_nu}  holds true, see Proposition 2.7 in \cite{BandiniRusso2}. On the other hand, assume that Hypothesis \ref{H_nu} is satisfied. Then, recalling Remark \ref{R_pred_supp}-(ii)-(iv)
		 and taking into account the additivity of the predictable projection operator, we have 
	\begin{align*}
	\one_{ K}= {}^{p}\left(\one_{D}\right) ={}^{p}\left(\one_{ K}\right) + {}^{p}\left(\one_{D \setminus K}\right) = \one_{ K} + {}^{p}\left(\one_{D \setminus K}\right)
	\end{align*}
	so that ${}^{p}\left(\one_{D \setminus K}\right)=0$. It follows that the predictable support of $D \setminus  K$ is an evanescent set. By Remark \ref{R_pred_supp}-(v) $D \setminus  K$ a thin set, therefore  $D \setminus  K= \cup_n [[T_n^i]]$, with $(T^i_n)_n$ disjoint totally inaccessible times, see  Corollary 5.43 in \cite{chineseBook} and  Remark at page 122 in \cite{chineseBook}.
\end{remark}

\begin{remark}\label{R_H_mu_BIS} If $\nu(\{t\}\times\R) =0$ for every $t \geq 0$,  then 	 
		 $J = K = \emptyset$ 
	and  
		Hypothesis 
		\ref{H_nu} trivially holds. 
\end{remark}

\begin{proposition}[Proposition 2.12 in \cite{BandiniRusso1}] 
\label{P_identiy_measures}  
	Let  $\mu$ be an integer-valued random measure and $X$  be an adapted c\`adl\`ag process, such that  
$(X,\mu)$  verifies Hypothesis \ref{H_X_mu}.
	Then, there exists a null set $\mathcal N$ such that, 
	for every Borel function $\varphi: [0,\,T] \times  \R \rightarrow \R_+$ satisfying  $\varphi(s,0)=0$, $s \in [0,\,T]$, we have, for every $\omega \notin \mathcal N$, 
	\begin{equation*}
	\int_{]0,\,\cdot]\times \R}\varphi(s,x)\,\mu^X(\omega,ds\,dx)= \int_{]0,\,\cdot]\times \R}\varphi(s,\tilde \gamma(\omega,s,e))\,\mu(\omega,ds\,de) + \sum_{0<s \leq \cdot} \varphi(s,\Delta X^p_s(\omega)).
	\end{equation*} 
	\end{proposition}
\begin{proposition}
\label{L_ident_mu_muX} 
Let  $\mu$ be an integer-valued random measure with compensator $\nu$ satisfying  Hypothesis \ref{H_nu}, and let $X$  be an adapted c\`adl\`ag process such that   and 
$(X,\mu)$  verifies Hypothesis \ref{H_X_mu}.
	Let  $\varphi: \Omega \times [0,\,T] \times \R \rightarrow \R_+$ 
	be a $\mathcal {\tilde P}$-measurable function such that $\varphi(\omega, s,0)=0$ for every $s \in [0,\,T]$, up to indistinguishability, and assume that there exists  a $\mathcal {\tilde P}$-measurable subset $A$ of $\Omega \times [0,\,T] \times \R$ satisfying 
	\begin{equation}\label{cond_int_A} 
	|\varphi|\,\one_{A} \star \mu^X \in \mathcal A^+_{\rm loc},\quad |\varphi|^2\,\one_{A^c} \star \mu^X \in \mathcal A^+_{\rm loc}. 
	\end{equation} 
	Then 
	\begin{equation}\label{Id_intstoch2} 
	\int_{]0,\,\cdot]\times \R}\varphi(s,x)\,\,(\mu^X-\nu^X)(ds\,dx) 
	= \int_{]0,\,\cdot]\times \R}\varphi(s,\tilde \gamma(s,e))\,(\mu-\nu)(ds\,de).
\end{equation}
\end{proposition}
\begin{remark}\label{R:migliorato}The result above consistently improves the achievements of Proposition 2.17  in  \cite{BandiniRusso2}. As a matter of fact,    condition \eqref{CondInMeno} 
	is no longer asked here.  This  
	 allows to  solve the identification problems under more general assumptions, see Theorems \ref{P_ident} and  \ref{P_ident_C0}, and therefore  it  extends the applicability of our results, see e.g. Section \ref{S:NewJumpDiff}.
\end{remark}
\proof	
Clearly the result holds if we show that $\varphi$ verifies \eqref{Id_intstoch2} 
under  one of the two following assumptions: 
	(i) $|\varphi| \star \mu^X \in \mathcal A_{\rm loc}^+$, 
	(ii) $|\varphi|^2 \star \mu^X \in \mathcal A_{\rm loc}^+$. 
By localization arguments, it is enough to show it when $|\varphi| \star \mu^X \in \mathcal A^+$, $|\varphi|^2 \star \mu^X \in \mathcal A^+$. 

\medskip 

\noindent \emph{Case  $|\varphi|  \star \mu^X \in \mathcal A^+ $.}
We will separate the proof into the following steps.
\begin{enumerate}
\item  Assume that $\theta: \Omega \times [0,\,T] \times \R \rightarrow \R_+$  
	be a $\mathcal {\tilde P}$-measurable function such that $\theta(\omega, s,0)=0$ for every $s \in [0,\,T]$,   	satifying \begin{equation}\label{Id_intstoch} 
	\int_{]0,\,\cdot]\times \R}\theta(s,x)\,\,(\mu^X-\nu^X)(ds\,dx) 
	= \int_{]0,\,\cdot]\times \R}\theta(s,\tilde \gamma(s,e))\,(\mu-\nu)(ds\,de) + \Gamma_\cdot^\theta 
	\end{equation} 
	with  $\Gamma^\theta$  a predictable  
	process. Then,  	 $\Gamma^\theta= 0$.  As a matter of fact,  $\Gamma^\theta$   is a predictable local martingale, and therefore continuous, see Remark 4 pag 194, Chapter VII,  in \cite{chineseBook}. On the other hand, being also a
	purely discontinuous martingale,  it follows that  it is the null process.

		\item Let $\phi: \Omega \times [0,\,T] \times \R \rightarrow \R_+$ 
	be a $\mathcal {\tilde P}$-measurable function such that $\phi(\omega, s,0)=0$ for every $s \in [0,\,T]$, and $|\phi|  \star \mu^X \in \mathcal A^+ $.
	Assume that $\phi= \phi \one_K$. Then \eqref{Id_intstoch2} holds true for $\varphi= \phi$. 
	\item Let $\psi: \Omega \times [0,\,T] \times \R \rightarrow \R_+$ 
	be a $\mathcal {\tilde P}$-measurable function such that $\psi(\omega, s,0)=0$ for every $s \in [0,\,T]$, and  $|\psi|  \star \mu^X \in \mathcal A^+ $.
	Assume that $\psi \one_K=0$. Then \eqref{Id_intstoch2} holds true for $\varphi= \psi$.
	\item Let $\varphi: \Omega \times [0,\,T] \times \R \rightarrow \R_+$ 
	be a $\mathcal {\tilde P}$-measurable function such that $\varphi(\omega, s,0)=0$ for every $s \in [0,\,T]$, and $|\varphi|  \star \mu^X \in \mathcal A^+ $.
Then $\varphi= \phi + \psi$, with  $\phi:= \varphi \one_K$,  $\psi = \varphi- \varphi \one_K$. 
In particular, $\phi=\phi \one_K$ and $\psi\one_K=0$, and $|\phi|  \star \mu^X \in \mathcal A^+ $,  $|\psi|  \star \mu^X \in \mathcal A^+ $ . By Steps 2. and 3. and the additivity property of the stochastic integral, it follows that \eqref{Id_intstoch2} holds true for $\varphi$. 
\end{enumerate}
It remains to prove Steps 2. and 3.

\noindent \emph{Step 2.} 
  We have   
\begin{align*}
	\int_{]0,\,\cdot]\times \R}\phi(s,x)\,\mu^X(ds\,dx)&=\int_{]0,\,\cdot]\times \R}\phi(s,x)\,\one_K(s)\,\mu^X(ds\,dx)\notag\\
	&=  \sum_{s \leq \cdot} \phi(s,\Delta X_s)\,\one_K(s)= \sum_{s \leq \cdot}\phi(s,\Delta X^i_s + \Delta X^p_s)\,\one_K(s),
	\end{align*}
where in the latter equality we have used that   $(X, \mu)$ satisfies  Hypotheses \ref{H_X_mu} with $X= X^i + X^p$. Since  $X^i$ is a c\`adl\`ag quasi-left-continuous process,   
$\Delta X^i_S=0$ for all predictable random
time $S$.
Recalling Remark \ref{R_pred_supp}-(iii) and Hypothesis \ref{H_nu}, we have 
\begin{align*}
\sum_{s \leq \cdot}\phi(s,\Delta X^i_s + \Delta X^p_s)\,\one_{K}(s)&=\sum_{s\leq \cdot}\phi(s,\Delta X^i_{s} + \Delta X^p_{s})\,\one_{\cup_n [[R_n]]}(s)\\
&=\sum_{n: \, R_n \leq \cdot}\phi(R_n,\Delta X^i_{R_n} + \Delta X^p_{R_n})\,\one_{[[R_n]]}(s)\\
&=\sum_{n: \, R_n \leq \cdot}\phi(R_n,\Delta X^p_{R_n})\,\one_{[[R_n]]}(s)=\sum_{s \leq \cdot}\phi(s,\Delta X^p_s)\,\one_{K}(s),
\end{align*}
so that \eqref{E_claim}  yields 
\begin{equation}\label{Ebis_claim} 
	\int_{]0,\,\cdot]\times \R}\phi(s,x)\,\mu^X(ds\,dx)=  \sum_{0<s \leq \cdot} \phi(s,\Delta X^p_s).
	\end{equation}
	By Proposition \ref{P_identiy_measures} together with \eqref{Ebis_claim}, 
	$
\int_{]0,\,\cdot]\times \R}\phi(s,\tilde \gamma(s,e))\,\mu(ds\,de)  =0
$. 
Therefore, being $\phi$ a non negative function,
$\int_{]0,\,\cdot]\times \R}\phi(s,\tilde \gamma(s,e))\,\nu(ds\,de)  =0
$. It follows that 
		$$
\int_{]0,\,\cdot]\times \R}\phi(s,\tilde \gamma(s,e))\,(\mu-\nu)(ds\,de)  =0.
$$
Adding and subtracting $\int_{]0,\,\cdot]\times \R}\phi(s,x)\,\nu^X(ds\,dx)$ in  \eqref{Ebis_claim} we get that \eqref{Id_intstoch}  holds for $\phi$ with 
$$
\Gamma^\phi := \sum_{0<s \leq \cdot} \phi(s,\Delta X^p_s) - \int_{]0,\,\cdot]\times \R}\phi(s,x)\,\nu^X(ds\,dx). 
$$
Being $\Gamma^\phi$ predictable, by Step 1 it follows that $\Gamma^\phi=0$, so that  
$$
\int_{]0,\,\cdot]\times \R}\phi(s,x)\,\,(\mu^X-\nu^X)(ds\,dx) =0.
$$
It follows that \eqref{Id_intstoch2} holds true for $\phi$, and reads $0=0$. 

\noindent \emph{Step 3.} Since $(X, \mu)$ satisfies  Hypotheses \ref{H_X_mu} and $J=K$ by  Hypothesis \ref{H_nu}, $\{\Delta X^p \neq 0 \} \subset K$,  so that $\one_{\{\Delta X^p \neq 0\}} \leq   \one_K$. We have 
$$
\sum_{0<s \leq \cdot} \psi(s,\Delta X^p_s)=\sum_{0<s \leq \cdot} \psi(s,\Delta X^p_s)\one_{\{\Delta X^p \neq 0\}}(s)\leq  \sum_{0<s \leq \cdot} \psi(s,\Delta X^p_s)\one_K(s)=0.
$$
Therefore, $\sum_{0<s \leq \cdot} \psi(s,\Delta X^p_s)=0$, and the equality in Proposition \ref{P_identiy_measures}	reads
\begin{equation}\label{E_claim} 
	\int_{]0,\,\cdot]\times \R}\psi(s,x)\,\mu^X(ds\,dx)= \int_{]0,\,\cdot]\times \R}\psi(s,\tilde \gamma(s,e))\,\mu(ds\,de).
	\end{equation}  
	Adding and subtracting $\int_{]0,\,\cdot]\times \R}\psi(s,x)\,\nu^X(ds\,dx)$ (resp. $\int_{]0,\,\cdot]\times \R}\psi(s,\tilde\gamma(s,e)\,\nu(ds\,de)$) in the left-hand side of \eqref{E_claim}  (resp. in the right-hand side of \eqref{E_claim}),  we get that \eqref{Id_intstoch}  holds for $\psi$ with 
$$
\Gamma^\psi := \int_{]0,\,\cdot]\times \R}\psi(s,\tilde\gamma(s,e)\,\nu(ds\,de)- \int_{]0,\,\cdot]\times \R}\psi(s,x)\,\nu^X(ds\,dx). 
$$
Being $\Gamma^\psi$ predictable, by Step 1 it follows that $\Gamma^\psi=0$, so that \eqref{Id_intstoch2} holds true for $\psi$. This concludes the proof in the case $|\varphi|  \star \mu^X \in \mathcal A^+$.

\medskip 

\noindent \emph{Case $|\varphi|^2  \star \mu^X \in \mathcal A^+ $.} This will follow from the previous one by approaching in $\mathcal L^2(\mu^X)$ the function $\varphi$ with $\varphi_{\varepsilon}(s,x):=\varphi(s,x)\,\one_{\varepsilon < |x| \leq 1/\varepsilon}\,\one_{s \in [0,\,T]}$. Indeed, $\varphi_\varepsilon(s,x) \star \mu^X \in \mathcal A^+$, by Cauchy-Schwarz inequality,
taking into account the fact that $\mu^X$, restricted to $\varepsilon \leq |x| \leq 1/\varepsilon$, is finite, since $\mu^X$ is $\sigma$-finite on $[0,\infty) \times \R$. The proof is done along the same steps as above, with the following slight modifications. 

\noindent \emph{Step 2'.} Set  $\phi_{\varepsilon}(s,x):=\phi(s,x)\,\one_{\varepsilon < |x| \leq 1/\varepsilon}\,\one_{s \in [0,\,T]}$, and notice that  $\phi_{\varepsilon}(s,x) = \phi_{\varepsilon}(s,x) \one_K$. Applying  Step 2 with  $\phi=\phi_{\varepsilon}$, we get that 
\begin{align}\label{PDmtg}
\int_{]0,\,\cdot]\times \R}\phi_{\varepsilon}(s,\tilde \gamma(s,e))\,(\mu-\nu)(ds\,de)  
=\int_{]0,\,\cdot]\times \R}\phi_{\varepsilon}(s,x)\,\,(\mu^X-\nu^X)(ds\,dx)=0.
\end{align}
We remind that, by \eqref{C(W)BIS}-\eqref{bracket_mart}-\eqref{G2norm},  if $||\phi_{\varepsilon}(s,x)-\phi(s,x)||^2_{\mathcal{G}^2(\mu^X)}$ and $||\phi_{\varepsilon}(s,\tilde \gamma(s,e))-\phi(s,\tilde \gamma(s,e))||^2_{\mathcal{G}^2(\mu)}$ converges to zero as $\varepsilon$ goes to zero, then \eqref{Id_intstoch2} holds for $\varphi$ replaced by $\phi$,  and reads $0=0$. 
Recalling Remark \ref{L_G2mu_L2mu_inclusion},  
we have 
	$||\phi_{\varepsilon}(s,x)-\phi(s,x)||^2_{\mathcal{G}^2(\mu^X)} \leq ||\phi_{\varepsilon}(s,x)-\phi(s,x)||^2_{\mathcal{L}^2(\mu^X)}$ and  $||\phi_{\varepsilon}(s,\tilde \gamma(s,e))-\phi(s,\tilde \gamma(s,e))||^2_{\mathcal{G}^2(\mu)} \leq ||\phi_{\varepsilon}(s,\tilde \gamma(s,e))-\phi(s,\tilde \gamma(s,e)))||^2_{\mathcal{L}^2(\mu)}$.  
By  the Lebesgue theorem, 
and  the fact that $\phi_{\varepsilon}$ converges pointwise to $\phi$,  
we have that  $||\phi_{\varepsilon}(s,x)-\phi(s,x)||^2_{\mathcal{L}^2(\mu^X)} \rightarrow 0$ and  $||\phi_{\varepsilon}(s,\tilde \gamma(s,e))-\phi(s,\tilde \gamma(s,e)))||^2_{\mathcal{L}^2(\mu)} \rightarrow 0$, and the conclusion follows. 

\noindent\emph{Step 3'.} Set  $\psi_{\varepsilon}(s,x):=\psi(s,x)\,\one_{\varepsilon < |x| \leq 1/\varepsilon}\,\one_{s \in [0,\,T]}$, and notice that  $\psi_{\varepsilon}(s,x) = \psi_{\varepsilon}(s,x) \one_K$. Arguing as in Step 3, we get that \eqref{Id_intstoch2} holds true for $\psi_{\varepsilon}$, namely
$$
	\int_{]0,\,\cdot]\times \R}\psi_{\varepsilon}(s,x)\,\,(\mu^X-\nu^X)(ds\,dx) 
	= \int_{]0,\,\cdot]\times \R}\psi_{\varepsilon}(s,\tilde \gamma(s,e))\,(\mu-\nu)(ds\,de).
$$
	The conclusion follows  
	arguing as in Step 2'. 
\endproof

\section{The identification problem}\label{Sec_application_BSDE}
In the present section we address the identification problem in two cases, the first one consisting in Theorem \ref{P_ident} and the second one consisting in Theorem \ref{P_ident_C0}.

Let  $(\Omega,\mathcal{F}, (\mathcal{F}_t)_{t \geq 0},\P)$ be a filtered probability space.
%
Let $\mu$ be an integer-valued random measure defined on $[0,T]\times \R$, with compensator $\nu$, 
and let $M$ be a 
local martingale, 
 with $M_0=0$. 
Let 
$\zeta$ be  a non-decreasing adapted 
c\`adl\`ag
process.
We will focus on 
 BSDEs driven by a compensated random measure $\mu-\nu$ of the form
  \eqref{GeneralBSDE}.
Here 
$\xi$ is an $\mathcal F_T$-measurable  square integrable random variable, 
$\tilde f: \Omega \times [0,\,T] \times \R^4 \rightarrow \R$
 is a measurable function, whose domain is equipped with the $\sigma$-field 
 $\mathcal F\otimes \mathcal B([0,\,T] \times \R^4)$.
A solution  of BSDE \eqref{GeneralBSDE} is a triple of processes $(Y,Z,U)$ such that  the first two  integrals in \eqref{GeneralBSDE} exist and are finite in the Lebesgue sense, $Y$ is adapted and c\`adl\`ag, $Z$ is  progressively measurable with  $Z \in L^2([0,\,T], d \langle M\rangle_t)$ a.s., and  
   $U \in \mathcal{G}^2_{\rm loc}(\mu)$.
   
\begin{remark}\label{R_uniq_G2}
Uniqueness 
means the following: 
if $(Y,Z,U)$, $(Y',Z',U')$ are solutions of the BSDE 
\eqref{GeneralBSDE}, then $Y=Y'$ in the sense of indistinguishability, 
$Z=Z'$ $d \P\,d \langle M\rangle_t$ a.e., and $U_t(e)-U'_t(e)$ in the sense of $ \mathcal{G}^2_{\rm loc}(\mu)$, namely, 
there is a predictable process $(l_t)$ such that 
$U_t(e)-U'_t(e) = l_t\,\one_K(t)$, $d\P\,\nu(dt\,de)$-a.e.	
The latter fact is a direct consequence of Proposition \ref{P_forBSDEs2}. 
In particular, if  $K = \emptyset$, 
then the third component  of the  BSDE solution is 
uniquely characterized in  
  $\mathcal{L}^2(\mu)$. 
\end{remark}

\begin{theorem}\label{P_ident} 
	Let $\mu$  be a random measure with compensator $\nu$
	satisfying Hypothesis \ref{H_nu}, and assume that 
	 $X$ is a c\`adl\`ag process such that $(X,\mu)$ verifies  Hypothesis \ref{H_X_mu}. 
	Let  $(Y,Z,U)$ be a solution to the BSDE \eqref{GeneralBSDE} such that the pair $(X, Y)$ satisfies 
	Hypothesis \ref{H_chain_rule_C01} with corresponding  function $v$.
Let $X^c$ denote the continuous local martingale of $X$ given in the canonical decomposition \eqref{Mc+Md+A}. 
	Then,  the pair $(Z,U)$ 
	fulfills 
	\begin{equation}\label{id_2BIS} 
	Z_t = \partial_x v(t,X_t) \,\frac{d\langle X^c,M\rangle_t}{d\langle M \rangle_t}\quad 
	d\P \,d\langle M \rangle_t \,\textup{-a.e.,} 
	\end{equation} 
	\begin{equation}\label{id_3} 
	\int_{]0,\,t]\times \R} H_s(e)\,(\mu-\nu)(ds\,de)=0,\quad  \forall \,\, t \in ]0,\,T],\,\, \textup{a.s.} 
	\end{equation} 
	with	 
	$H_s(e):= U_s(e)-(v(s,X_{s-}+ \tilde\gamma(s,e))-v(s,X_{s-}))$. 
	
\noindent 	If, in addition, $H \in \mathcal G^2_{\rm loc}(\mu)$, 
	then there exists a predictable process $(l_s)$ such that 
\begin{equation}\label{Hl_id}
	H_s(e)= l_s\,\one_{K}(s),\quad d\P\,\nu(ds\,de)\textup{-a.e.}
\end{equation}
\end{theorem} 
\begin{remark}\label{R:id_part}
In particular, it follows from \eqref{Hl_id} and Hypothesis \ref{H_nu} for $\mu$,  that  
$H_s(e)= 0$, $d\P\,\nu^c(ds\,de)\textup{-a.e.}$
and	
$H_s(e)= l_s$, 
$d\P\,\nu^d(ds\,de)\textup{-a.e.}$
\end{remark}

\proof
By assumption, the couple $(X,Y)$ satisfies Hypothesis \ref{H_chain_rule_C01} with corresponding function $v$.
We are thus in the condition to apply 
Theorem \ref{T_C1_dec_specialweak_Dir} 
to $v(t,\,X_t)$. 
We set 
$\varphi(s,x):=v(s,X_{s-}+x)-v(s,X_{s-})$. 
Since $X$ is of finite quadratic variation and verifies \eqref{Sec:WD_CNS}, and $X$ and $v$ satisfy {\eqref{Sec:WD_E_C}, by 
\eqref{A2_Aloc_NEW} and \eqref{F2_Aloc_NEW}
 we see that the process $\varphi$ verifies condition \eqref{cond_int_A} with $A=\{|x| >1\}$. 
Moreover $\varphi(s,0)=0$.   Since  $\mu$ verifies Hypothesis \ref{H_nu} and  $(X,\mu)$ verifies Hypothesis \ref{H_X_mu}, we can apply Proposition \ref{L_ident_mu_muX} to $\varphi(s,x)$. Identity \eqref{C01_special_WD_formula} in Theorem  \ref{T_C1_dec_specialweak_Dir}   becomes 
\begin{align} 
v(t,\,X_t) 
&=v(0,X_0) + \int_{]0,\,t]\times \R} (v(s,X_{s-}+ \tilde \gamma(s,e))-v(s,X_{s-}))\,(\mu-\nu)(ds\,de) 
\nonumber\\ 
&\,\,+ \int_{]0,\,t]}\partial_x v(s,X_{s})\,d X^{c}_s+ A^v(t), \label{dec_Y} 
\end{align} 
 where $A^v$ is a predictable $(\mathcal F_t)$-orthogonal  process. In particular, $v(t, X_t)$ is a special weak Dirichlet process. 
On the other hand, the process $Y_t = v(t,X_t)$ fulfills the BSDE \eqref{GeneralBSDE}. In particular it is a special semimartingale, and therefore a special weak Dirichlet process. 
By Remark 
\ref{Sec:WD_P_unique_decomp}, which states the uniqueness of the decomposition of a special weak Dirichlet process, 
 we get \eqref{id_3}
and 
\begin{equation}\label{Id2} 
\int_{]0,\,t]} Z_s \, dM_s=\int_{]0,\,t]} \partial_x v(s,X_{s})\,d X_s^c. 
\end{equation} 
In particular, from \eqref{Id2} we get 
\begin{align*} 
0 &=\langle \int_{]0,\,t]} Z_s dM_s-\int_{]0,\,t]} \partial_x v(s,X_{s})\,d X_s^c,\, M_t \rangle\\ 
&=\int_{]0,\,t]} Z_s d\langle M \rangle_s -\int_{]0,\,t]} \partial_x v(s,X_{s})\,\frac{d \langle X^c,\, M \rangle_s}{d \langle M\rangle_s}\,d \langle M\rangle_s\\ 
&=\int_{]0,\,t]} \bigg(Z_s  - \partial_x v(s,X_{s})\,\frac{d \langle  X^c,\, M \rangle_s}{d \langle M\rangle_s}\bigg)\,d \langle M\rangle_s, 
\end{align*} 
that gives identification \eqref{id_2BIS}. 
If in addition $H\in \mathcal G^2_{\rm loc}(\mu)$, the predictable bracket at time $t$ of the purely discontinuous martingale in   identity \eqref{id_3} is well-defined, and 
   by \eqref{bracket_mart}  equals $C(H)$ given in \eqref{C(W)BIS}. 
Since $C(H)_T=0$ a.s., the conclusion follows from Proposition
\ref{P_forBSDEs2}.
\endproof
If $\mu= \mu^X$, Theorem \ref{P_ident}  simplifies in the following way. \begin{theorem}\label{C_ident_1} 
	Let $X$ be a  c\`adl\`ag  process, whose   jump measure $\mu^X$ 
	with compensator $\nu^X$
	satisfies Hypothesis \ref{H_nu}. 
	Let  $(Y,Z,U)$ be a solution to the BSDE \eqref{GeneralBSDE} with $\mu=\mu^X$,  such that the pair $(X, Y)$ satisfies 
	Hypothesis \ref{H_chain_rule_C01} with corresponding  function $v$.
Let $X^c$ denote the continuous local martingale of $X$ given in the canonical decomposition \eqref{Mc+Md+A}. 
	 	Then,  the pair $(Z,U)$ 
	fulfills 
	\begin{equation}\label{id_2} 
	Z_t = \partial_x v(t,X_t) \,\frac{d\langle X^c,M\rangle_t}{d\langle M \rangle_t}\quad 
	d\P \,d\langle M \rangle_t \,\textup{-a.e.,} 
	\end{equation} 
	\begin{equation}\label{id_3BIS} 
	\int_{]0,\,t]\times \R} H_s(x)\,(\mu^X-\nu^X)(ds\,dx)=0,\quad  \forall \,\, t \in ]0,\,T],\,\, \textup{a.s.} 
	\end{equation} 
	with	 
	$H_s(x):= U_s(x)-(v(s,X_{s-}+x)-v(s,X_{s-}))$. 
	
\noindent 	If, in addition, $H \in \mathcal G^2_{\rm loc}(\mu^X)$, 
	then there exists a predictable process $(l_s)$ such that 
\begin{equation}\label{Hl_idBIS}
	H_s(x)= l_s\,\one_{K}(s),\quad d\P\,\nu^X(ds\,dx)\textup{-a.e.}
\end{equation}
\end{theorem} 
\proof
The proof goes along the same lines of the one of Theorem  \ref{P_ident}, the only difference being  that we replace \eqref{dec_Y} directly with identity \eqref{C01_special_WD_formula} in Theorem  \ref{T_C1_dec_specialweak_Dir}.
\endproof

Let us now consider a BSDE  driven only by a purely discontinuous martingale,  of the form
\begin{align}\label{GeneralBSDE_disc} 
Y_t &= \xi + \int_{]t,\,T]\times \R} \tilde f(s,\,e,\,Y_{s-},\,U_{s}(e))\,  d\zeta_s - \int_{]t,\,T]\times \R} U_s(e)\,(\mu-\nu)(ds\,de). 
\end{align} 
\begin{theorem}\label{P_ident_C0} 
	Let  $\mu$ be a random measure with compensator $\nu$ satisfying Hypothesis \ref{H_nu}, and assume that
	  $X$ is a process such that $(X,\mu)$ verifies Hypothesis \ref{H_X_mu}. 
	Let   $(Y,U)$ be a solution to the BSDE \eqref{GeneralBSDE_disc}, 
	such that $(X,Y)$ satisfies Hypothesis  \ref{H_chain_rule_C0} with corresponding function $v$.
	Then,  the random field  $U$ 	satisfies 
	\begin{equation}\label{id_3disc} 
	\int_{]0,\,t]\times \R} H_s(e)\,(\mu-\nu)(ds\,de)=0 \quad \forall t \in ]0,\,T],\,\,\textup{a.s.} 
	\end{equation} 
	with 
	$H_s(e):=U_s(e)-(v(s,X_{s-}+ \tilde\gamma(s,e))-v(s,X_{s-}))$.

\noindent 	If, in addition, $H \in \mathcal G^2_{\rm loc}(\mu)$, 
	then there exists a predictable process $(l_s)$ such that 
\begin{equation}\label{id_l_2}
	H_s(e)= l_s\,\one_{K}(s),\quad d\P\,\nu(ds\,de)\textup{-a.e.}
\end{equation}
\end{theorem}

\proof
Set 
$ 
\varphi(s,x):=v(s,X_{s-}+x)-v(s,X_{s-})$. 
By condition (ii) in Hypothesis \ref{H_chain_rule_C0}, the process $\varphi$ verifies condition \eqref{cond_int_A} with $A=\Omega \times [0,\,T] \times \R$. 
Moreover $\varphi(s,0)=0$.   Since  $\mu$ verifies Hypothesis \ref{H_nu}, 
and  $(X, \mu)$ verifies Hypothesis \ref{H_X_mu}, 
we can apply Proposition \ref{L_ident_mu_muX} to $\varphi(s,x)$. Identity \eqref{C00_decomp} in Proposition \ref{P_C00_chain_rule}  becomes 
\begin{align} 
v(t,\,X_t) 
=v(0,X_0) + \int_{]0,\,t]\times \R} (v(s,X_{s-}+ \tilde \gamma(s,e))-v(s,X_{s-}))\,(\mu-\nu)(ds\,de) 
+ A^v(t),\label{dec_YC0} 
\end{align}
 where $A^v$ is a predictable $(\mathcal F_t)$-orthogonal  process.  
At this point we recall that the process $Y_t = v(t,X_t)$ fulfills BSDE \eqref{GeneralBSDE_disc}. 
Again, the uniqueness of   a special weak Dirichlet process (see Remark \ref{Sec:WD_P_unique_decomp}) 
yields   identity 
\eqref{id_3disc}.  
If in addition we assume that $H\in \mathcal G^2_{\rm loc}(\mu)$, the predictable bracket at time $t$ of the purely discontinuous martingale in   identity \eqref{id_3disc} is well-defined,  and 
   by \eqref{bracket_mart} equals $C(H)$ given in  \eqref{C(W)BIS}. 
Since $C(H)_T=0$ a.s.,
the conclusion follows from Proposition
\ref{P_forBSDEs2}. 
\endproof
Also in this case, the result simplifies when $\mu= \mu^X$. 
\begin{theorem}\label{C_ident_C0} 
Let $X$ be a  c\`adl\`ag  process, whose   jump measure $\mu^X$ 
	with compensator $\nu^X$
	satisfies Hypothesis \ref{H_nu}.
	Let   $(Y,U)$ be a solution to the BSDE \eqref{GeneralBSDE_disc} with $\mu=\mu^X$, 
	such that $(X,Y)$ satisfies Hypothesis  \ref{H_chain_rule_C0} with corresponding function $v$. 
	Then,  the random field  $U$ 	satisfies 
	\begin{equation}\label{id_3discBIS} 
	\int_{]0,\,t]\times \R} H_s(x)\,(\mu^X-\nu^X)(ds\,dx)=0 \quad \forall t \in ]0,\,T],\,\,\textup{a.s.} 
	\end{equation} 
	with 
	$H_s(x):=U_s(x)-(v(s,X_{s-}+ x)-v(s,X_{s-}))$.

\noindent 	If, in addition, $H \in \mathcal G^2_{\rm loc}(\mu^X)$, 
	then there exists a predictable process $(l_s)$ such that 
\begin{equation}\label{id_l_2BIS}
	H_s(x)= l_s\,\one_{K}(s),\quad d\P\,\nu^X(ds\,dx)\textup{-a.e.}
\end{equation}
\end{theorem} 
\proof
The proof goes along the same lines of the one of Theorem  \ref{P_ident_C0}, the only difference being  that we replace \eqref{dec_YC0} directly with identity \eqref{C01_special_WD_formula} in Theorem  \ref{T_C1_dec_specialweak_Dir}.
\endproof

\section{Applications}

 \subsection{The jump diffusion case}\label{S:NewJumpDiff}
 
Let  $(\Omega,\mathcal{F}, (\mathcal{F}_t)_{t \geq 0},\P)$ be a filtered probability space.
In the present section we consider a random measure $\mu$ and a process $X$ satisfying the following. 
\begin{hypothesis}\label{HXMU}
$\mu(ds\,de)$  is an integer-valued random measure with  compensator $\nu(ds\,de)=d A_s \,\phi_s(de)$  satisfying Hypothesis \ref{H_nu}. 
$X$  is a solution of the equation 
	\begin{equation}\label{X_SDE} 
	X_t = X_0 +\int_0^t b(s,X_{s-})\,d C_s +\int_0^t \sigma(s,X_s)\,dN_s + \int_{]0,\,t]\times \R} \gamma(s,X_{s-},e)\,(\mu-\nu)(ds\,de). 
	\end{equation} 
Here	$N$ is a continuous martingale,  $C$ is an increasing predictable c\`adl\`ag  process, with $C_0=0$, such that 
	$\{\Delta C \neq 0\} 
	\subset J$.  
Moreover 	 $b,\sigma: \Omega \times [0,\,T]\times \R \rightarrow \R$, $\gamma:\Omega \times [0,\,T]\times \R \times \R \rightarrow \R$ are  $\tilde {\mathcal P}$-measurable maps such that 
	 $(\omega,s,e)\mapsto \gamma(\omega, s,X_{s-}(\omega),e)\in \mathcal G^1_{\rm loc}(\mu)$ and 
	\begin{align} 
	&\int_0^t |b(s,X_{s-})|\, d C_s < \infty \,\,\textup{a.s.},\label{Ver_1}\\
	&\int_0^t |\sigma(s,X_s)|^2\,d \langle N\rangle_s < \infty\,\, \textup{a.s.},
	\label{Ver_2}\\ 
	&\gamma(s, X_{s-}, e)\one_K(s) \equiv 0.\label{Ver_3} 
	\end{align}	 
\end{hypothesis}
 We have the following results. 
\begin{lemma}\label{Ex_guida}
Let  $X$ be a c\`adl\`ag process and  $\mu$ be a random measure such that $(X,\mu)$  satisfies  Hypothesis \ref{HXMU}. 
Then  $(X,\mu)$  satisfies Hypothesis \ref{H_X_mu} with decomposition $X=X^i+X^p$, where
	\begin{align} 
	X^i_t &= \int_{]0,\,t]\times \R} \tilde  \gamma(s, e)
	\,(\mu-\nu)(ds\,de),\label{Xi}\\ 
	X^p_t &= X_0 +\int_0^t b(s,X_{s-})\,d C_s +\int_0^t \sigma(s,X_s)\,dN_s
	,\label{Xp} 
	\end{align} 
	with $\tilde \gamma(\omega, s,e) = \gamma(\omega, s,X_{s-}(\omega),e)
	$.

\end{lemma} 
 \begin{remark}
 This result extends  Lemma 2.19 in \cite{BandiniRusso2} in  two ways. Firstly, we allow the coefficients $b$, $\gamma$ and $\sigma$ to be random.  Secondly, we no longer ask condition \eqref{CondInMeno} on $\mu$, and instead we  ask condition \eqref{Ver_3} on the coefficient $\gamma$ in \eqref{X_SDE}. This allows for instance  to consider the case when $\mu$ does not fulfills condition \eqref{CondInMeno} and $\gamma$ satisfies \eqref{Ver_3}, which  was not included in Lemma 2.19 in \cite{BandiniRusso2}. 
 \end{remark}
 \proof
By  \eqref{X_SDE}-\eqref{Xi}-\eqref{Xp}, together with \eqref{Ver_3},  it straightly follows that $X= X^i+ X^p$. 
Let us now show that 
 $X^i$ and $X^p$  in \eqref{Xi}-\eqref{Xp} are respectively  a c\`adl\`ag  quasi-left-continuous  adapted process and a c\`adl\`ag   predictable process. 	 
The fact that  $X^p$ is predictable straight follow from \eqref{Xp}. Concerning $X^i$,   it is enough to prove that  $\Delta X^i_S\,\one_{\{S< \infty\}} = 0$ a.s., for any  $S$  predictable time, see Remark \ref{R:totinac_leftcont}. 
Recalling that by Hypothesis \ref{H_nu} we have $J=K$, and that $\gamma$ fulfills \eqref{Ver_3}, we get 
\begin{equation} \label{gammatilde}
\Delta X^i_s = \int_{\R} 
\,\gamma(s,X_{s-},e)\,\one_{D \setminus K}(s) \,
\mu(\{s\},de). 
\end{equation} 
 Recalling \eqref{mubeta}, 
\eqref{gammatilde} can be rewritten as 
\begin{align}
\Delta X^i_s(\omega) &=\gamma(\omega, s,X_{s-}(\omega),\beta_s(\omega))\,\,\one_{D \setminus K}(\omega, s),\label{DeltaXiDmenoK}\\
&=  \gamma(\omega, s,X_{s-}(\omega),\beta_s(\omega))\,\,\one_{\cup_n [[T_n^i]]}(\omega, s),\label{DeltaXi}
\end{align}
where the second line follows by  the fact that 
$D\setminus K = \cup_n [[T_n^i]]$ up to an evanescent set, $(T_n^i)_n$ being a sequence of totally inaccessible times with disjoint graphs, see   Remark \ref{R_Hp_3.1}.
Identity \eqref{DeltaXi}  gives, for any $S$ finite predictable time,  
\begin{align*}
\Delta X^i_S(\omega)\,\one_{\{S< \infty\}} 
&=\gamma(\omega,S,X_{S-}(\omega),\beta_S(\omega))\, \sum_n \one_{[[T_n^i]]}(\omega, S(\omega))\,\one_{\{S< \infty\}}
\end{align*}
which is zero 
being $(T_n^i)_n$ a sequence  of totally inaccessible times. 
The fact that $\{\Delta X^i \neq 0\} \subset D$ also directly follows from \eqref{gammatilde}. 
To prove that $\Delta X^i_s(\omega) = \tilde \gamma(\omega, s,\cdot)$, $dM^{\P}_{\mu}(\omega,s)$-a.e., it is enough to show that 
$$ 
\sper{\int_{]0,\,T]\times \R}\mu(\omega, ds\,de)\,|\tilde \gamma(\omega, s,e)-\Delta X^i_s(\omega)|\,\one_{\tilde \Omega_n}(\omega,s)}=0. 
$$ 
By the structure of $\mu$ it follows that, for every $n \in \N$, 
\begin{align*} 
\sper{\int_{]0,\,T]\times \R}\mu(\omega, ds\,de)\,|\tilde \gamma(\omega, s,e)-\Delta X^i_s(\omega)|\one_{\tilde \Omega_n}(\omega,s)} \leqslant
\sum_{s \in ]0,\,T]} \sper{\one_{D}(\cdot,s)\,|\tilde \gamma(\cdot, s,\beta_s(\cdot))-\Delta X^i_s(\cdot)|}, 
\end{align*} 
which vanishes taking into account \eqref{DeltaXiDmenoK} and that $\gamma$ fulfills \eqref{Ver_3}. 
Finally, 
since $N$ is continuous, it follows from \eqref{Xp} 
that 
\begin{equation}\label{DeltaXp} 
\Delta X^p_s = b(s,X_{s-})\,\Delta C_s, 
\end{equation} 
 so that 
$ 
\{\Delta X^p \neq 0\}\subset \{\Delta C \neq 0\} 
\subset J
$.
\endproof 

\begin{lemma}\label{L_BBP_EC} 
	Let  $X$ be a c\`adl\`ag process and  $\mu$ be a random measure such that $(X,\mu)$  satisfies  Hypothesis \ref{HXMU}. 
	Assume 
 that 
\begin{equation}\label{sum_Aloc2} 
\sum_{s \in]0,\,\cdot]}|b(s, X_{s-})|^2 |\Delta C_s|^2  + \int_{]0,\cdot]\times \R}|\gamma(X_{s-},e)|^2\,\nu(ds\,de)\in \mathcal A_{\rm loc}^+. 
\end{equation} 	
	Then the following holds.
	\begin{itemize}
		\item [(i)]
$X$ is  a special weak Dirichlet process with  finite quadratic variation 
such that $(X,\mu)$ verifies condition  \eqref{Sec:WD_CNS}; 
	\item[(ii)] 	if  
	$v:[0,\,T] \times \R \rightarrow \R$ is a function of  $C^{0,1}$ class such that $x \mapsto \partial_x v(s,x)$ has linear growth, uniformly in $s$,  
condition \eqref{Sec:WD_E_C} holds for $X$ and $v$. 
	\end{itemize} 
\end{lemma} 
\proof 
(i) By \eqref{X_SDE}-\eqref{Ver_1}-\eqref{Ver_2},  $X$ is a special semimartingale.
   In particular condition \eqref{Sec:WD_CNS} holds, see  Corollary 11.26 in \cite{chineseBook}. Moreover, obviously 
  $X$ has finite quadratic variation.
\noindent (ii) For some constant $c$ we have
\begin{align} 
&\int_{]0,\cdot]\times \R} |v(s,X_{s-} +x )-v(s,X_{s-})-x\,\partial_x v(s,X_{s-})|\,\one_{\{|x| >1\}}\,\mu^X(ds\,dx)\nonumber\\ 
&=\sum_{0<s \leq \cdot} |v(s,X_{s})-v(s,X_{s-})-\partial_x v(s,X_{s-})\,\Delta X_s|\,\one_{\{|\Delta X_s| >1\}}\nonumber\\ 
&\leq \sum_{0<s \leq \cdot}|\Delta X_s|\,\one_{\{|\Delta X_s| >1\}}\,\left(\int_0^1|\partial_{x} v(s,X_{s-}+ a\,\Delta X_s)|\,da +\int_0^1 |\partial_{x} v(s,X_{s-})|\,da\right)\nonumber \\ 
&\leq 
2\,c\,\int_{]0,\cdot]\times \R}|X_{s-}|\,|x|\,\one_{\{|x| >1\}} \,\mu^X(ds\,dx) + \sum_{s \leq \cdot}|\Delta X_s|^2\,\one_{\{|\Delta X_s| >1\}}.\label{verificaAloc} 
\end{align} 
The first term in the right-hand side of \eqref{verificaAloc} belongs to $\mathcal A^+_{\rm loc}$, 
taking into account \eqref{Sec:WD_CNS} and the fact  that $X_{s-}$ is locally bounded 
being c\`agl\`ad. 

On the other hand, since $X$ is of finite quadratic variation, by Lemma 2.10-(ii) in \cite{BandiniRusso1} we have that $\sum_{s \in]0,\,T]}|\Delta X_s|^2 < \infty$, a.s.
Consequently, the second term in the right-hand side of \eqref{verificaAloc} belongs to $\mathcal A_{\rm loc}^+$ 
if we prove that
\begin{equation}\label{sum_Aloc} 
\sum_{s \in]0,\,\cdot]}|\Delta X_s|^2\in \mathcal A_{\rm loc}^+.
\end{equation} 
By \eqref{DeltaXi}-\eqref{DeltaXp},   $\Delta X_s= b(s, X_{s-}) \Delta C_s + \int_{\R}\gamma(X_{s-},e)\,\mu(\{s\},\,de)$, so that
\begin{align*}
\sum_{s \in]0,\,\cdot]}|\Delta X_s|^2 
&\leq   \sum_{s \in]0,\,\cdot]}|b(s, X_{s-})|^2 |\Delta C_s|^2  +  \int_{]0,\cdot]\times \R}|\gamma(X_{s-},e)|^2\,\mu(ds\,de), 
\end{align*}
which belongs to  $\mathcal A_{\rm loc}^+$ because of 
\eqref{sum_Aloc2}.
\endproof 
Let $W$ be a Brownian motion and $\mu(ds\,de)$ be  a
random measure with compensator $\nu(ds\,de)=\phi_s(de) dA_t$.
We will  focus on 
	the BSDE  
	\begin{align}\label{BSDE_BaBuPa} 
	Y_t &= g(X_T) + \int_{]t,\,T]} f(s,\,X_s,\,Y_s,\,Z_s,\,U_{s}(\cdot))\,  d A_s- \int_{]t,\,T]} Z_s \, d W_s - \int_{]t,\,T]\times \R} U_s(e)\,(\mu-\nu)(ds\,de), 
	\end{align} 
	which constitutes	a particular case of the BSDE  \eqref{GeneralBSDE}. 
	The process $X$ appearing in \eqref{BSDE_BaBuPa}   is a solution to \eqref{X_SDE}  satisfying \eqref{Ver_1}-
\eqref{Ver_2}-
\eqref{Ver_3}
-\eqref{sum_Aloc2}.
BSDEs of the type  \eqref{BSDE_BaBuPa} are  considered in \cite{papapantoleon_possamai_saplaouras};  under suitable assumptions, the existence and uniqueness of a solution $(Y,Z,U)\in \mathcal{S}^2 \times \mathcal{L}^2\times \mathcal{G}^2(\mu)$ is  established.

We are ready to give the identification result  in the present framework. 
\begin{corollary}\label{C_id_BBP} 
	Let   $(Y,Z,U)\in 	\mathcal{S}^2 \times \mathcal{L}^2\times \mathcal{G}^2(\mu)$ be 
	a solution  to the BSDE \eqref{BSDE_BaBuPa}. 
	Then  the pair $(Z,U)$ 
	satisfies 
	\begin{align}
& Z_t = \sigma(X_t)\,\partial_x u(t,X_t) \quad d\P\, dt\textup{-a.e.,} \label{Z_BBPnew} \\
	& \int_{]0,\,t]\times \R} H_s(e)\,(\mu-\nu)(ds\,de)=0, \quad \forall t \in ]0,\,T],\,\, \textup{a.s.}, \label{id_3BBPnew} 
	\end{align} 
	where 
	$
	H_s(e):=U_s(e)-(v(s,X_{s-}+ \gamma(s,X_{s-},e)
	)-v(s,X_{s-})). 
	$

\noindent 	If in addition $H \in \mathcal G^2_{\rm loc}(\mu)$,  
	\begin{align} 
U_s(e)&=v(s,X_{s-}+ \gamma(s,X_{s-},e))-v(s,X_{s-}), \quad d\P\,\nu^c(ds\,de)\textup{-a.e.}\label{U_BBPnew}
	\end{align}
and there exists a predictable process $(l_s)$ such that 
	\begin{align} 
U_s(e)&= l_s,
 \quad d\P\,\nu^d(ds\,de)\textup{-a.e.}\label{U_BBPnew2} 
	\end{align}	 
\end{corollary}	 
\proof 
By Lemma \ref{L_BBP_EC}, 
the pair $(X, Y)$ verifies Hypothesis \ref{H_chain_rule_C01}. 
On the other hand, by Lemma \ref{Ex_guida},    $(X,\mu)$ verifies  Hypothesis \ref{H_X_mu} with $X^i$ and $X^p$ given respectively by \eqref{Xi} and \eqref{Xp}.
We can then apply Theorem \ref{P_ident}: since $X^c_\cdot=\int_0^\cdot \sigma(X_t)\,dW_t$ and $M=W$, 
\eqref{id_2BIS}    gives \eqref{Z_BBPnew}, while \eqref{id_3}
 with $\tilde \gamma(s,e) = \gamma(s,X_{s-}(e))
 $
 yields \eqref{id_3BBPnew}. 
If in addition  $H \in \mathcal{G}^2(\mu)$,    
\eqref{U_BBPnew}-\eqref{U_BBPnew2} follows by \eqref{Hl_id} and Remark \ref{R:id_part}. 
\endproof 
\begin{remark}
This result significantly extends Corollary 4.3 in \cite{BandiniRusso2}, where $\mu(ds\,de)$ where a Poisson random measure with deterministic compensator $\nu(de)ds$ (and in particular condition \eqref{CondInMeno} held   being $K= \emptyset$). 	Here we insist on the fact that we  deal with a general random measure $\mu$ not necessarily verifying condition \eqref{CondInMeno}.
\end{remark}

\subsection{The PDMPs case}\label{S:The PDMPs case}
	We assume that $X$ is a  piecewise deterministic Markov process (PDMP)
generated by a marked point process $(T_n, \zeta_n)$, where $(T_n)_n$ are increasing random times such that 
	$ 
	T_n \in ]0,\,\infty[, 
	$ 
	where either there is a finite number of 
 times  $(T_n)_n$  or $\lim_{n \rightarrow \infty} T_n = + \infty$, 
	and $\zeta_n$ are random variables in $[0,1]$.
	We will follow the notations in \cite{Da-bo}, Chapter 2, Sections 24 and 26. 	The behavior of the PDMP $X$ is described by a triplet of 
	local characteristics  $(h,\lambda,Q)$: 
	$h: ]0,\,1[ \rightarrow \R$  is a Lipschitz continuous function, 
	$\lambda: ]0,1[ \rightarrow \R$ is a measurable function such that 
	$\sup_{x \in ]0,1[}|\lambda(x)| < \infty$, 
	and $Q$ is a transition probability  measure on $[0,1]\times\mathcal{B}(]-1,1[)$. 	Some other technical assumptions  are specified in the over-mentioned reference, that we do not recall here.
	Let us denote by $\Phi(s,x)$  the unique solution of $
	{g}'(s) = h(g(s))$, $g(0)= x$.
	The process $X$ can  be   defined as 
		\begin{equation} \label{X_eq}
			X(t)= 
	\left\{ 
	\begin{array}{ll} 
	\Phi(t,x),\quad t \in [0,\,T_{1}[\\ 
	\Phi(t-T_n, \zeta_n),\quad t \in [T_n,\,T_{n+1}[, 
	\end{array} 
	\right. 
	\end{equation}  
	and verifies the equation 
	\begin{equation}\label{PDP_dynamic_estended} 
	X_t=X_0 + \int_{0}^t h(X_s)\,ds + \int_{]0,\,t]\times \R}x\,\mu^X(ds\,dx)
	\end{equation} 
	with 
\begin{equation}\label{muPDMPS}
	\mu^X(ds\,dx) 
	=\sum_{n \geq 1} \one_{\{\zeta_{n}\in ]0,1[\}} \delta_{( T_n,\,\zeta_n- \zeta_{n-1})}(ds\,dx).  
\end{equation}
The knowledge of $(b,\,\lambda,\,Q)$ completely specifies the law of $X$, see Section 24 in \cite{Da-bo}, and also Proposition 2.1 in \cite{BandiniPDMPs}. In particular, let $\P$ be the unique probability measure under which  
	the compensator of $\mu^X$ has the form 
	\begin{equation}\label{nuPDPs} 
	\nu^X(ds\, dx) = (\lambda(X_{s-})\,ds + d p^{\ast}_s)\,Q(X_{s-},\,dx), 
	\end{equation} 
	where $\lambda$ has been trivially extended to $[0,1]$ by the zero value, and  
	\begin{equation}\label{p_ast} 
	p^{\ast}_t = \sum_{n =1}^{\infty} \one_{\{t \geq T_n\}}\,\one_{\{X_{{T_n}-} \in \{0,1\}\}}
	\end{equation} 
	is the predictable process counting the number of jumps of  $X$ from the  boundary of its domain. 
	 
	From \eqref{nuPDPs}, we can write decomposition $\nu^X(ds\, dx)= \phi_s(dx) dA_s$ 
	with $d A_s = \lambda(X_{s-})\,ds + d p^{\ast}_s$ and $\phi_{s}(dx)= Q(X_{s-},dx)$. 
	In particular, $A$ is predictable (not deterministic) and discontinuous, with jumps 
	$\Delta A_s(\omega) =
	 \Delta p^{\ast}_s(\omega)=\one_{\{X_{s-}(\omega) \in \{0,1\}\}}$. 
	Consequently, 
	\begin{equation}\label{PDP_J=K} 
	J=
	K =\{(\omega,t):X_{t-}(\omega)\in \{0,1\}\}.
	\end{equation}

\begin{remark} 
In \cite{BandiniRusso2}
  we  asked  the measure $\mu^X$ to satisfy condition \eqref{CondInMeno}, 
   with  entails  
the existence of a  function $\beta: \{0,1\}\rightarrow ]-1,1[$ such that 
	\begin{equation*}
	Q(y,\,dx)=\delta_{\beta(y)}(dx)\quad \textup{a.s.}, 
	\end{equation*} 
	see Lemma 4.11 in \cite{BandiniRusso2}.  In the present paper  we can 
	avoid this assumption  and  work
  with the whole class of PDMPs.
\end{remark}

	Let us  consider a BSDE driven by the compensated random measure $\mu^X-\nu^X$, where $\mu^X$ is the integer-valued random measure  in   \eqref{muPDMPS}  associated to a piecewise deterministic Markov process $X$   with values in the interval $[0,1]$,
	of the   form 
	\begin{equation}\label{BSDE_PDP} 
	Y_t = g(X_T) + \int_{]t,\,T]} f(s,\,X_{s-},\,Y_{s-},\,U_{s}(\cdot))\,  d A_s - \int_{]t,\,T]\times \R} U_s(e)\,(\mu^X-\nu^X)(ds\,de). 
	\end{equation} 
Existence and uniqueness results for  solutions $(Y,U)\in 	\mathcal{L}^2 \times  \mathcal{G}^2(\mu)$ to  BSDEs driven by  purely discontinuous martingales (that include \eqref{BSDE_PDP}  as a special case) were  established under suitable assumptions in \cite{BandiniBSDE} and in the recent work \cite{BandiniSICON}.

\begin{lemma}\label{L_Bandini_EC0} 
We set $E = [0,\,1]$.
	Let   $Y$ be a  special semimartingale such that its  martingale component is purely discontinuous. 	
Let $X$	be a  c\`adl\`ag process with values in $E$, with a  finite number of jumps on each compact interval. 
	Assume that $Y_t = v(t,X_t)$ for some  
	 function   $v:[0,\,T] \times \R \rightarrow \R$ such that its restriction to $[0,\,T]\times E$ is continuous. 
	Then $(X,Y)$ satisfies Hypothesis \ref{H_chain_rule_C0} with corresponding function $v$. 
\end{lemma} 
\proof	 
The proof is  the same as the one of Lemma 4.13 in \cite{BandiniRusso2}. 
\endproof

\begin{corollary}\label{C_id_PDPs} 
	Let  $X$  be a PDMP with jump measure $\mu^X$ with compensator $\nu^X$ given by \eqref{nuPDPs},  
	and let  $(Y,U)\in 	\mathcal{L}^2 \times  \mathcal{G}^2(\mu^X)$ be a solution  to the BSDE \eqref{BSDE_PDP}.
	Assume that $Y_t = v(t,X_t)$ for some continuous function  $v$. 
	Then  the random field $U$		 satisfies
	\begin{equation}\label{id_PDMP} 
	\int_{]0,\,t]\times \R} H_s(x)\,(\mu^X-\nu^X)(ds\,dx)=0, \quad \forall \,t \in ]0,\,T],\,\, \textup{a.s.} 
	\end{equation} 
with $H_s(x) := U_s(x)-(v(s,X_{s-} + x)-v(s,X_{s-}))$.	
	If in addition 
	$H_s(x)\in \mathcal G^2_{\rm loc}(\mu^X)$, 
	\begin{equation}\label{final_ident} 
	U_s(x)=v(s,X_{s-} + x)-v(s,X_{s-}) \quad d\P\,\lambda(X_{s-})\,Q(X_{s-},\,dx)\,\one_{X_{s-} \in ]0,\,1[} ds\textup{-a.e.} 
	\end{equation} 
and  there exists a predictable process $(l_s)$ such that 
\begin{equation}\label{final_ident_nud}
U_s(x)= v(s,X_{s-} + x)-v(s,X_{s-}) +l_s
, \quad d\P\,
Q(X_{s-},\,dx)\one_{X_{s-} \in \{0,\,1\}}\,dp^\ast_s\textup{-a.e.}
\end{equation}
\end{corollary} 
\proof 
Hypothesis \ref{H_chain_rule_C0} holds for $(X,Y)$ by Lemma \ref{L_Bandini_EC0}. 
We are in condition to apply Theorem \ref{C_ident_C0}, which gives   \eqref{id_PDMP}.
 If, in addition, $H \in \mathcal G^2_{\rm loc}(\mu^X)$, by \eqref{id_l_2BIS} together with Remark \ref{R:id_part} with $\mu=\mu^X$,
$$
H_s(x)= 0, \quad d\P\,\nu^{X,c}(ds\,de)\textup{-a.e.},\\
$$
and	there exists a predictable process $(l_s)$ such that 
$$
H_s(x)= l_s, \quad d\P\,\nu^{X,d}(ds\,de)\textup{-a.e.}
$$ 
At this point, recalling \eqref{nuPDPs}, and being   $J=K$, we see that 
\begin{align*}
\nu^{X,c}(ds\,dx)
&=\nu^X(ds\,dx)\,\one_{K^c}(s)=\lambda(X_{s-})\, Q(X_{s-},dx)\,\one_{K^c}(s)\,ds, \\
{\nu}^{X,d}(ds \,dx)&=\nu^X(ds\,dx)\,\one_{K}(s)=Q(X_{s-},\,dx)\,\one_K(s) \,dp^\ast_s. 
\end{align*}
Then, since by \eqref{PDP_J=K}  we have $K =\{(\omega,s): X_{s-}(\omega) \in \{0,1\}\}$,   and \eqref{final_ident}-\eqref{final_ident_nud} follow.
\endproof 

\subsection{The  jump-diffusion case  with distributional drift}\label{S:The  jump-diffusion case  with distributional drift}

Let $\sigma, \beta \in C^0$ such that $\sigma >0$. We consider formally a PDE 
operator, obtained  by mollification  (see e.g. \cite{frw1}, \cite{frw2}),   
of the type
\begin{equation}\label{Lbeta}
 L \psi = \frac{1}{2}\sigma^2  \psi'' + \beta' \psi'. 
\end{equation}
  \begin{hypothesis}\label
  {H:h2}
 We assume the existence of a function 
 $\Sigma(x) := \lim_{n \rightarrow \infty} 2 \int_0^x \frac{\beta'_n}{\sigma_n^2} (y) dy$
in  $C^0$, independently from the mollifier.
Moreover $\Sigma \in C^{\alpha}$ for some $\alpha \in [0,\,1]$,  the function $\Sigma$ is lower bounded,
and  
 $$
 \int_{-\infty}^0 e^{-\Sigma(x)} dx = \int_{0}^{+\infty} e^{-\Sigma(x)} dx= + \infty.
 $$
 \end{hypothesis}

 \begin{definition}
 We will denote by 	$\mathcal D_{L}$ the set of all $f \in C^1$ such that there exists some $\dot{l} \in C^0$ with $L f = \dot{l}$ in the sense of \cite{frw1}.
  This defines without ambiguity $L : \mathcal D_{L} \rightarrow C^0$.
 \end{definition}
We introduce the following definition of martingale problem. For an increasing process $A$, we will use the notation $A_t = A_t ^c + \sum_{s \leq t} \Delta A_s$. We also  denote by $\mathcal C^+(\R)$ the set of bounded Borel functions of $\R$, vanishing inside  a neighborhood of $0$. 
In particular, if any two positive measures $\eta$, $\eta'$ on $\R$ with $\eta(\{0\})=\eta'(\{0\})=0$, and $\eta(x: |x|>\varepsilon) <\infty$, $\eta'(x: |x|>\varepsilon) <\infty$ are such that $\eta(f) = \eta'(f)$ for all $f \in \mathcal C^+(\R)$, then $\eta= \eta'$.
 \begin{definition}\label{D:new_mrtg_pb}
 A couple $(X, \P)$ is said to solve the martingale problem 
related to a given  operator ${L}$ of the form \eqref{Lbeta}, a random measure $\nu(ds\,dx) = \phi_s(dx)\, dA_s $,  an initial condition $X_0=x_0\in \R$,  and a  domain $\mathcal { D}\subset \mathcal {D}_L$,  if
the following holds.
 \begin{itemize}
  	\item[(i)] 
for any  $f \in \mathcal {D}$,  
\begin{equation}\label{int_cond_fX}
\int_0^t\int_{\R} [f(X_{s-} + x) -f(X_{s-})-x\,f'(X_{s-})
]\,\nu(ds\,dx) \in \mathcal A_{\textup{loc}}^+,
\end{equation}  
with respect to  $\P$;
\item [(ii)] for any  $f \in \mathcal {D}$, 
the process
\begin{equation*}
Z^f := f(X_{\cdot}) - f(X_0) - \int_0^{\cdot} {L} f(X_s) dA_s^c
-  \int_{]0,\cdot]}\int_{\R} [f(X_{s-} + x) -f(X_{s-})-x\,f'(X_{s-})\one_{J^c}]\,\nu(ds\,dx)
\end{equation*}
is a local martingale  under $\P$;
\item[(iii)]  for every $g \in \mathcal C^+(\R)$,  
$G^g:=g \ast \mu^X -g \ast \nu$
  is a local martingale under $\P$.
  \end{itemize}
\end{definition}

\begin{remark}\label{R:D}
Assume that $(\nu, \P)$ satisfies 
condition \eqref{alpha+1integ} for som $\alpha \in [0,\,1]$.  Then condition \eqref{int_cond_fX} holds with respect to  $\P$ for any $f \in \mathcal{ D}$ with 
\begin{equation}\label{barD}
	\mathcal{D}:=\{f \in 
\mathcal{D}_{L}
 \,\,\textup{with} \,\,f \in C^{1+ \alpha}_{\textup loc},  f'\,\, \textup{bounded}\}.
\end{equation}
\end{remark}

Let $L$ be an operator of the form \eqref{Lbeta}, for which  Hypotheses \ref{H:h2} holds, and  
let $\nu(ds\,dx) = \phi_s(dx) \, d A_s$ be a predictable random measure.
Let $(X, \P)$ be   a  solution to the martingale problem in  Definition \ref{D:new_mrtg_pb} related to $X_0$, $L$,   
 $\nu(ds\,dx)$, and $\mathcal{ D}$ given in \eqref{barD}, and such  that      $(\nu, \P)$ satisfies 
 condition \eqref{alpha+1integ}.
 
 The following result is given 
 in \cite{BandiniRusso_DistributionalDrift}, where we study under suitable assumptions the well-posedness of the martingale problem in Definition \ref{D:new_mrtg_pb}. 
 \begin{proposition}\label{P:DRIFTDISTRIBUZIONALE}  
 	$X$ is a special weak Dirichlet process
(with respect to its canonical filtration) 
 of finite quadratic variation with canonical decomposition $X= X_0 + M^{X}  + \Gamma$, with $\Gamma$  a predictable and  $\mathcal F_t^X$-orthogonal process, and $M^{X}=  M^{X,d} + X^c$, satisfying condition \eqref{Sec:WD_CNS} under $\P$. 
In particular,  $\nu$ is the $\P$-compensator of $\mu^X$, and 
\begin{align*}
	M^{X,d} = \int_{]0,\,\cdot]}\int_{\R} x \,(\mu^X-\nu)(ds\,dx),\quad 
	\langle X^c\rangle =\int_0^{\cdot}\sigma^2(X_s)\,d A^c_s. 
\end{align*}
\end{proposition}
We are interested in   BSDEs under $\P$  driven by the compensated random 	measure $\mu^X - \nu$ and  the continuous martingale $X^c$, of the   form 
	\begin{align}\label{BSDE_distrdrift} 
	Y_t &= g(X_T) + \int_{]t,\,T]} f(s,\,X_{s},Y_s, \,Z_{s},\,U_{s}(\cdot))\,  d A_s \notag\\
	&
	- \int_{]t,\,T]} Z_s \frac{1}{\sigma(X_s)}d X_s^c
	- \int_{]t,\,T]\times \R} U_s(x)\,(\mu^X-\nu)(ds\,dx).
	\end{align} 
	A consequence of our identification Theorem  \ref{C_ident_1} is the following.
\begin{corollary}\label{C_id_distr} 
Assume that $\mu^X(ds\,dx)$
satisfies Hypothesis \ref{H_nu}.
  Let   $(Y,Z,U)\in 	\mathcal{S}^2 \times \mathcal{L}^2\times \mathcal{G}^2(\mu^X)$ be a solution  to the BSDE \eqref{BSDE_distrdrift} 
Assume that   $Y_t = v(t,X_t)$ for some deterministic function $v: [0,T] \times \R \rightarrow \R$ of class $C^{0,1}$ such that $v$ and $X$ verify   condition \eqref{Sec:WD_E_C} under $\P$. 
	Then  the pair $(Z,U)$ 
	satisfies 
	\begin{equation}\label{Z_BBP} 
Z_t = \sigma(X_t)\,\partial_x v(t,X_t) \quad d\P\, d A^c_t\textup{-a.e.,} 
	\end{equation} 
	\begin{equation}\label{id_3BBP} 
	\int_{]0,\,t]\times \R} H_s(x)\,(\mu^X-\nu)(ds\,dx)=0, \quad \forall \,t \in ]0,\,T],\,\, \textup{a.s.} 
	\end{equation} 
	with
	$H_s(x):=U_s(x)-(v(s,X_{s-}+ x
	)-v(s,X_{s-}))$. 
If, in addition, $H \in \mathcal G^2_{\rm loc}(\mu^X)$, 
	then there exists a predictable process $(l_s)$ such that 
\begin{equation}\label{Hl_idDIST}
	H_s(x)= l_s\,\one_{K}(s),\quad d\P\,\nu(ds\,dx)\textup{-a.e.}
\end{equation}	
\end{corollary}	 
\proof 
We aim at applying Theorem \ref{C_ident_1}. 
By assumption $\mu^X$  satisfies Hypothesis \ref{H_nu}.
On the other hand, by Proposition \ref{P:DRIFTDISTRIBUZIONALE}, 	$X$ is a  special weak Dirichlet process      of finite quadratic variation  with its canonical decomposition  $X=X^c+M^d+\Gamma$.
 In addition,  by assumption $X$ satisfies condition \eqref{Sec:WD_CNS} under $\P$, and  condition \eqref{Sec:WD_E_C} holds for $X$ and  $v$ under $\P$,  see  Remark \ref{R:D}. This implies the validity of Hypothesis \ref{H_chain_rule_C01} for $(X,Y)$. 
We can then apply Theorem \ref{C_ident_1}: since $\langle X^c_\cdot\rangle=\int_0^\cdot \sigma^2(X_t)\,d A^c_t$ and $M= \int_0^\cdot \frac{1}{\sigma(X_s)}d X_s^c$, 
formula \eqref{id_2}   gives \eqref{Z_BBP}, while \eqref{id_3BIS}
yields \eqref{id_3BBP}.
If in addition  $H \in \mathcal{G}^2(\mu^X)$,   
then \eqref{Hl_idDIST} 
 follows by \eqref{Hl_idBIS}, recalling that $\nu$ is the $\P$ compensator of $\mu^X$. 
\endproof

\small
\paragraph{Acknowledgements.} 
The two authors benefited of the support of  GNAMPA  
  project  \emph{Controllo ottimo stocastico con osservazione parziale: metodo di randomizzazione ed equazioni di 
Hamilton-Jacobi-Bellman sullo spazio di Wasserstein}. The work of the first named author was partially  supported by PRIN 2015 
\emph{Deterministic and Stochastic Evolution equations}. 
The work of the second named author was partially supported by a public grant as part of the
{\it Investissement d'avenir project, reference ANR-11-LABX-0056-LMH,
  LabEx LMH,}
in a joint call with Gaspard Monge Program for optimization, operations research and their interactions with data sciences.
 
\addcontentsline{toc}{chapter}{Bibliography} 
\bibliographystyle{plain} 
\bibliography{BiblioLivreFRPV_TESI} 

\def\cprime{$'$} \def\cprime{$'$} \def\cprime{$'$} \def\cprime{$'$}
\begin{thebibliography}{10}

\bibitem{BandiniBSDE}
\textsc{Bandini, E.}
\newblock Existence and uniqueness for backward stochastic differential
  equations driven by a random measure.
\newblock {\em Electronic Communications in Probability}, 20(71):1--13, 2015.

\bibitem{BandiniPDMPs}
\textsc{Bandini, E.}
\newblock Optimal control of {P}iecewise-{D}eterministic {M}arkov {P}rocesses:
  a {BSDE} representation of the value function.
\newblock {\em ESAIM: Control, Optimization and Calculus of Variations},
  4(1):311--354, 2018.

\bibitem{BandiniSICON}
\textsc{Bandini, E.}
\newblock Constrained {BSDE}s driven by a non quasi-left-continuous random
  measure and optimal control of {PDMP}s on bounded domains.
\newblock {\em SIAM Journal on Control and Optimization}, 57(6):3767--3798,
  2019.

\bibitem{BandiniConfortola}
\textsc{Bandini, E. \& Confortola, F.}
\newblock Optimal control of semi-{M}arkov processes with a backward stochastic
  differential equations approach.
\newblock {\em Mathematics of Control, Signals and Systems}, 29(1):1--35, 2017.

\bibitem{BandiniFuhrman}
\textsc{Bandini, E. \& Fuhrman, M.}
\newblock Constrained {BSDE}s representation of the value function in optimal
  control of pure jump {M}arkov processes.
\newblock {\em Stochastic Process. Appl.}, 127(5):1441--1474, 2017.

\bibitem{BandiniRusso_DistributionalDrift}
\textsc{Bandini, E. \& Russo, F.}
\newblock The martingale problem related to an {SDE} with jumps and
  distributional drift.
\newblock {\em Work in progress}.

\bibitem{BandiniRusso1}
\textsc{Bandini, E. \& Russo, F.}
\newblock Weak {D}irichlet processes with jumps.
\newblock {\em Stochastic Processes and their Applications},
  127(12):4139--4189, 2017.

\bibitem{BandiniRusso2}
\textsc{Bandini, E. \& Russo, F.}
\newblock Special weak {D}irichlet processes and {BSDE}s driven by a random
  measure.
\newblock {\em Bernoulli}, 24(4A):2569--2609, 2018.

\bibitem{BaBuPa}
\textsc{Barles, G. \& Buckdahn, R. \& Pardoux, E.}
\newblock Backward stochastic differential equations and integral-partial
  differential equations.
\newblock {\em Stochastic Rep.}, 60(1-2):57--83, 1997.

\bibitem{barrassopreprint1}
\textsc{Barrasso, A. \& Russo, F.}
\newblock Backward {S}tochastic {D}ifferential {E}quations with no driving
  martingale, {M}arkov processes and associated {P}seudo {P}artial
  {D}ifferential {E}quations.
\newblock 2017.
\newblock Preprint, hal-01431559, v1.

\bibitem{barrasso3}
\textsc{Barrasso, A. \& Russo, F.}
\newblock Martingale driven {BSDE}s, {PDE}s and other related deterministic
  problems.
\newblock 2017.
\newblock Preprint, hal-01566883.

\bibitem{BP}
\textsc{Buckdahn, R. \& Pardoux, E.}
\newblock {BSDE}'s with jumps and associated integral-stochastic differential
  equations.
\newblock {\em Preprint}, 1994.

\bibitem{CarboneFerrSantacroce}
\textsc{Carbone, R. \& Ferrario, B. \& Santacroce, M.}
\newblock Backward stochastic differential equations driven by c\`adl\`ag
  martingales.
\newblock {\em Theory Probab. Appl.}, 52:304--314, 2008.

\bibitem{CoFu-m}
\textsc{Confortola, F. \& Fuhrman, M.}
\newblock Backward stochastic differential equations associated to {M}arkov
  jump processes and applications.
\newblock {\em Stochastic Processes and their Applications}, 124:289--316,
  2014.

\bibitem{cjms}
\textsc{Coquet, F. \& Jakubowski, A. \& M{\'e}min, J. \& S{\l}omi{\'n}ski, L.}
\newblock Natural decomposition of processes and weak {D}irichlet processes.
\newblock In {\em In memoriam {P}aul-{A}ndr\'e {M}eyer: {S}\'eminaire de
  {P}robabilit\'es {XXXIX}}, volume 1874 of {\em Lecture Notes in Math.}, pages
  81--116. Springer, Berlin, 2006.

\bibitem{Da-bo}
\textsc{Davis, M. H. A.}
\newblock {\em Markov models and optimization.}, volume~49 of {\em Monographs
  on Statistics and Applied Probability}.
\newblock Chapman $\&$ Hall., 1993.

\bibitem{er2}
\textsc{Errami, M. \& Russo, F.}
\newblock {$n$}-covariation, generalized {D}irichlet processes and calculus
  with respect to finite cubic variation processes.
\newblock {\em Stochastic Process. Appl.}, 104(2):259--299, 2003.

\bibitem{frw1}
\textsc{Flandoli, F. \& Russo, F. \& Wolf, J.}
\newblock Some {SDE}s with distributional drift. {I}. {G}eneral calculus.
\newblock {\em Osaka J. Math.}, 40(2):493--542, 2003.

\bibitem{frw2}
\textsc{Flandoli, F. \& Russo, F. \& Wolf, J.}
\newblock Some {SDE}s with distributional drift. {II}. {L}yons-{Z}heng
  structure, {I}t\^o's formula and semimartingale characterization.
\newblock {\em Random Oper. Stochastic Equations}, 12(2):145--184, 2004.

\bibitem{FuhrmanTessitore}
\textsc{Fuhrman, M. \& Tessitore, G.}
\newblock Generalized directional gradients, backward stochastic differential
  equations and mild solutions of semilinear parabolic equations.
\newblock {\em Appl Math Optim}, 108:263--298, 2003.

\bibitem{gr1}
\textsc{Gozzi, F. \& Russo, F.}
\newblock Verification theorems for stochastic optimal control problems via a
  time dependent {F}ukushima-{D}irichlet decomposition.
\newblock {\em Stochastic Process. Appl.}, 116(11):1530--1562, 2006.

\bibitem{gr}
\textsc{Gozzi, F. \& Russo, F.}
\newblock Weak {D}irichlet processes with a stochastic control perspective.
\newblock {\em Stochastic Process. Appl.}, 116(11):1563--1583, 2006.

\bibitem{chineseBook}
\textsc{He, S. \& Wang, J. \& Yan, J.}
\newblock {\em Semimartingale theory and stochastic calculus}.
\newblock Science Press Bejiing New York, 1992.

\bibitem{jacod_book}
\textsc{Jacod, J.}
\newblock {\em Calcul {S}tochastique et {P}robl\`emes de martingales}, volume
  714 of {\em Lecture Notes in Mathematics}.
\newblock Springer, Berlin, 1979.

\bibitem{JacodBook}
\textsc{Jacod, J. \& Shiryaev, A. N.}
\newblock {\em Limit theorems for stochastic processes}, volume 288 of {\em
  Grundlehren der Mathematischen Wissenschaften [Fundamental Principles of
  Mathematical Sciences]}.
\newblock Springer-Verlag, Berlin, second edition, 2003.

\bibitem{LaachirRusso}
\textsc{Laachir, I. \& Russo, F.}
\newblock B{SDE}s, c\`adl\`ag martingale problems, and orthogonalization under
  basis risk.
\newblock {\em SIAM J. Financial Math.}, 7(1):308--356, 2016.

\bibitem{papapantoleon_possamai_saplaouras}
\textsc{Papapantoleon A. \& Possamai D. \& Saplaouras A.}
\newblock Existence and uniqueness results for {BSDE} with jumps: the whole
  nine yards.
\newblock {\em Electron. J. Probab.}, 23(121):1--68, 2018.

\bibitem{pardouxpeng}
\textsc{Pardoux, {\'E}. \& Peng, S.}
\newblock Adapted solution of a backward stochastic differential equation.
\newblock {\em Systems Control Lett.}, 14(1):55--61, 1990.

\bibitem{PardouxPeng92}
\textsc{Pardoux, \'E. \& Peng, S.}
\newblock Backward {S}tochastic {D}ifferential {E}quations and {Q}uasilinear
  {P}arabolic {P}artial {D}ifferential {E}quations.
\newblock {\em Lecture Notes in CIS}, 176:200--217, 1992.

\bibitem{Peng91}
\textsc{Peng, S.}
\newblock {P}robabilistic {I}nterpretation for {S}ystems of {Q}uasilinear
  {P}arabolic {P}artial {D}ifferential {E}quations.
\newblock {\em Stochastics}, 37:61--74, 1991.

\bibitem{Peng92b}
\textsc{Peng, S.}
\newblock A generalized dynamic programming principle and
  {H}amilton-{J}acobi-{B}ellman {E}quation.
\newblock {\em Stochastics}, 38:119--134, 1992.

\bibitem{rv95}
\textsc{Russo, F. \& Vallois, P.}
\newblock The generalized covariation process and {I}t\^o formula.
\newblock {\em Stochastic Process. Appl.}, 59(1):81--104, 1995.

\bibitem{TaLi}
\textsc{Tang, S. J. \& Li, X. J.}
\newblock Necessary conditions for optimal control of stochastic systems with
  random jumps.
\newblock {\em SIAM J. Control Optim.}, 32:1447--1475, 1994.

\bibitem{xia}
\textsc{Xia, J.}
\newblock Backward stochastic differential equations with random measures.
\newblock {\em Acta Mathematicae Applicatae Sinica}, 16(3):225--234, 2000.

\end{thebibliography}

\end{document}